\documentclass[reqno,11pt]{amsart}

\usepackage{amssymb}
\usepackage[all,cmtip]{xy}
\usepackage[vcentermath]{youngtab}
\newtheorem{Proposition}{Proposition}[section]
\newtheorem{Definition}[Proposition]{Definition}
\newtheorem{Lemma}[Proposition]{Lemma}
\newtheorem{Theorem}[Proposition]{Theorem}
\newtheorem{Corollary}[Proposition]{Corollary}

\DeclareMathOperator{\Val}{Val}
\DeclareMathOperator{\vol}{vol}

\DeclareMathOperator{\im}{Im}
\DeclareMathOperator{\kl}{Kl}
\DeclareMathOperator{\Sym}{Sym}

\DeclareMathOperator{\ch}{char}
\DeclareMathOperator{\Gr}{Gr}
\DeclareMathOperator{\SL}{SL}
\DeclareMathOperator{\GL}{GL}
\DeclareMathOperator{\Spin}{Spin}
\DeclareMathOperator{\G}{G_2}
\DeclareMathOperator{\SU}{SU}
\DeclareMathOperator{\SO}{SO}
\DeclareMathOperator{\U}{U}
\DeclareMathOperator{\Sp}{Sp}
\DeclareMathOperator{\nc}{nc}
\newcommand{\R}{\mathbb{R}}
\newcommand{\C}{\mathbb{C}}
\newcommand{\K}{\mathcal{K}}
\newcommand{\F}{\mathbb{F}}
\newcommand{\h}{\mathbb{H}}
\newcommand{\spn}{\Sp(n)}
\newcommand{\spnu}{\Sp(n) \cdot \U(1)}
\newcommand{\spnsp}{\Sp(n) \cdot \Sp(1)}

\title{Invariant valuations on quaternionic vector spaces}
\author{Andreas Bernig}

\email{bernig@math.uni-frankfurt.de}

\address{Institut f\"ur Mathematik, Goethe-Universit\"at Frankfurt,
Robert-Mayer-Str. 10, 60054 Frankfurt, Germany}

\begin{document}

\begin{abstract}
The spaces of $\spn$-, $\spnu$- and $\spnsp$- invariant, translation invariant,
continuous convex valuations on the quaternionic vector space $\h^n$ are
studied. Combinatorial dimension formulas involving Young diagrams and Schur
polynomials are proved.  
\end{abstract}

\thanks{{\it MSC classification}:  53C65,  
52A22 
\\Supported by the Swiss National Science Foundation grants PP002-114715/1 and SNF
200020-121506/1.}

\maketitle 

\section{Introduction and statement of main results}
\label{sec_introduction}

One of the most important formulas in integral geometry is the principal
kinematic formula of Chern-Blaschke-Santal\'o. Let $V$ be an $m$-dimensional Euclidean vector space, let $\overline{\SO(V)}$ be the group generated by translations and rotations, endowed with an appropriate
Haar measure, and let $\omega_k$ being the volume of a $k$-dimensional unit ball. 

Let $K,L$ be compact convex sets. Then the principal
kinematic formula reads
\begin{equation} \label{eq_kin_formula}
\int_{\overline{\SO(V)}} \chi(K \cap \bar g L)d\bar g= \sum_{k=0}^m
\binom{m}{k}^{-1}\frac{\omega_k \omega_{m-k}}{\omega_m}\mu_k(K)\mu_{m-k}(L),
\end{equation}
where $\mu_0,\ldots,\mu_m$ are the {\it intrinsic volumes} (see \cite{klain_rota,schneider_book93}) and $\chi=\mu_0$ denotes the {\it Euler characteristic}, which is defined by $\chi(K)=1$ for all non-empty compact convex sets $K$ and $\chi(\emptyset)=0$. 

A nice conceptual proof due to Hadwiger goes as follows (see \cite{klain_rota} for details). First one notes that
the intrinsic volumes are {\it valuations}, i.e. finitely 
additive maps on the space of compact convex bodies. If $\Val^{\SO(V)}$ denotes 
the space of continuous, translation and rotation invariant valuations, then the
left hand side of formula \eqref{eq_kin_formula} with 
$K$ (or $L$) fixed belongs to $\Val^{\SO(V)}$. Using {\it Hadwiger's
theorem}, which states that $\mu_0,\ldots,\mu_m$ is a basis of $\Val^{\SO(V)}$, one obtains a formula of type \eqref{eq_kin_formula}, but with unknown
coefficients on the right hand side. These coefficients 
may be easily computed by the {\it template method}, i.e. by plugging in spheres on both sides of the formula and comparing coefficients.  

Instead of taking the full rotation group, one may restrict to some subgroup
$G$ of $\SO(V)$. Again we let $\bar G$ be the group generated by $G$ and translations. Since the main ingredient in the above proof 
is Hadwiger's theorem, an analogous {\it $G$-kinematic formula} exists for all
groups $G$ such that the space $\Val^G$ of $\bar G$-invariant continuous convex valuations is finite-dimensional. To write these formulas explicitly may be a challenge, since the template method is in general too weak to determine the constants. 

Alesker showed in \cite{alesker_survey07} that $\Val^G$ is
finite-dimensional if and only if $G$ acts transitively on the unit
sphere. Connected compact groups with this property were
classified by Borel and Montgomery-Samelson. There are six infinite
series  
\begin{equation} \label{eq_series}
 \SO(n), \U(n), \SU(n), \Sp(n), \Sp(n) \cdot \U(1), \Sp(n) \cdot \Sp(1) 
\end{equation}
and three exceptional groups 
\begin{equation} \label{eq_exceptions}
 \G, \Spin(7), \Spin(9). 
\end{equation}  

A valuation $\mu \in \Val^G$ is called homogeneous of degree $k$ if
$\mu(tK)=t^k\mu(K)$ for all compact convex sets and all $t \geq 0$. The
corresponding subspace is denoted by $\Val_k^G$. 

The classical, and most important, case is $G=\SO(n)$. By Hadwiger's theorem, we
have $\dim \Val_k^{\SO(n)}=1$ for $k=0,1,\ldots,n$. 

The case $G=\U(n)$ has been extensively studied in the last few years
\cite{tasaki03,alesker_mcullenconj01,
alesker03_un, fu06, bernig_fu_hig, abardia_thesis, abardia_gallego_solanes}. In
order to compare with our results in the quaternionic cases, we mention just
some results. 

Alesker showed in \cite{alesker_mcullenconj01} that 
\begin{displaymath}
 \dim \Val_k^{\U(n)}=\min\left\{\left\lfloor \frac{k}{2}\right\rfloor,
\left\lfloor \frac{2n-k}{2}\right\rfloor\right\}+1.
\end{displaymath}

Let us fix a sequence of holomorphic isometric embeddings $\C \hookrightarrow
\C^2 \hookrightarrow \C^3 \hookrightarrow \ldots$, they induce a
sequence of restriction maps 
\begin{displaymath}
\Val^{\U(1)} \leftarrow \Val^{\U(2)} \leftarrow \Val^{\U(3)} \leftarrow
\ldots
\end{displaymath}
Let 
\begin{displaymath}
 \Val^{\U(\infty)}:=\lim_{\longleftarrow} \Val^{\U(n)}
\end{displaymath}
be the inverse limit, which is a graded vector space (in fact a graded algebra). Elements of $\Val^{\U(\infty)}$ are called {\it global valuations}.

Alesker's results imply in particular that the restriction map 
\begin{equation} \label{eq_global_restriction_un}
 \Val_k^{\U(\infty)} \to \Val_k^{\U(n)}
\end{equation}
is surjective for all $n$ and $k$ and injective for $n \geq k$. Its kernel
was described by Fu \cite{fu06}. 

With respect to the product structure
introduced by Alesker
\cite{alesker04_product}, there is an isomorphism of graded algebras 
\begin{displaymath}
 \Val^{\U(\infty)} \cong \C[t,s],
\end{displaymath}
where $t,s$ are variables of degree $1$ and $2$ respectively.  In particular,
the Poincar\'e series of
$\Val^{\U(\infty)}$ is given by 
\begin{equation} \label{eq_poincare_series_un}
 \sum_{k=0}^\infty \dim \Val_k^{\U(\infty)}x^k= \frac{1}{(1-x)(1-x^2)}.
\end{equation}

Using the $\U(n)$-case, the integral geometry of the groups $\SU(n)$ was studied
in \cite{bernig_sun09}, and two of the three exceptional cases were
treated in
\cite{bernig_g2}. 

In the quaternionic cases (i.e. in the cases $G=\spn, \spnu, \spnsp$) relatively
little is known. The case $n=1$ was studied in \cite{alesker_su2_04} and
\cite{bernig_quat09}. Since $\Sp(1) \cong \SU(2)$, this is also a special case of
the general theory for $\SU(n)$, which was developed in \cite{bernig_sun09}.
Using plurisubharmonic functions in quaternionic variables, Alesker constructed
in \cite{alesker05} some $\spnsp$-invariant valuations on $\h^n$. For
instance, the {\it
Alesker-Kazarnovskii-pseudo-volume} (which was called {\it quaternionic
pseudo-volume} in \cite{alesker05}) is an element of $\Val_n^{\spnsp}$. 

Let us now describe our results in the quaternionic cases. Let $V \cong \h^n$ be
an $n$-dimensional quaternionic vector space. We consider $V$ as a {\it right}
vector space. The {\it compact symplectic group} $\Sp(n)$ then acts from the left
on $V$ by usual matrix multiplication. Moreover, the group $\Sp(1)$ of
quaternions of norm $1$ and its subgroup of complex numbers of norm $1$ act by
diagonal multiplication from the right on $V$. The subgroup of $\SO(4n)$
generated by $\spn$ and $\Sp(1)$ (resp. $\U(1)$) is denoted by $\spnsp$ (resp.
$\spnu$). More details on quaternionic vector spaces will be
given in
Section \ref{sec_back_rep}.

Our first main theorem is similar to \eqref{eq_global_restriction_un}, but the
proof will be different. Let $\h^1 \hookrightarrow \h^2 \hookrightarrow \h^3
\hookrightarrow \ldots$ be a sequence of quaternionic isometric embeddings. The corresponding restrictions yield a sequence
\begin{displaymath}
 \Val^{\Sp(1)} \leftarrow \Val^{\Sp(2)} \leftarrow \Val^{\Sp(3)} \leftarrow
\ldots.
\end{displaymath}
whose inverse limit is denoted by   
\begin{displaymath}
 \Val^{\Sp(\infty)}:=\lim_{\longleftarrow} \Val^{\Sp(n)}.
\end{displaymath}
The grading on each $\Val^{\Sp(n)}$ induces a grading on $\Val^{\Sp(\infty)}$. 

The spaces $\Val^{\Sp(\infty)\cdot \U(1)}$ and $\Val^{\Sp(n) \cdot \Sp(1)}$ are
defined in an analogous way. Elements in these spaces are called {\it global
valuations}, as opposed to {\it local valuations} in the spaces
$\Val^{\spn},\Val^{\spnu},\Val^{\spnsp}$. 

Our first main theorem is the quaternionic analogue of
\eqref{eq_global_restriction_un}. 

\begin{Theorem} \label{thm_local_global}
 The restriction maps 
\begin{align*}
 \Val_k^{\Sp(\infty)} & \to  \Val_k^{\Sp(n)}\\
 \Val_k^{\Sp(\infty)\cdot \U(1)} & \to \Val_k^{\Sp(n)\cdot \U(1)}\\
 \Val_k^{\Sp(\infty)\cdot \Sp(1)} & \to \Val_k^{\Sp(n) \cdot \Sp(1)}
\end{align*}
are surjective for all $k$ and $n$ and injective for $n \geq k$. 
\end{Theorem}

The second main theorem describes the dimensions of the spaces of global
valuations. Together with the previous theorem, it also yields the dimensions of
the spaces $\Val_k^{\spn}, \Val_k^{\spnu}, \Val_k^{\spnsp}$ when $k \leq n$
(and when $k \geq 3n$, see Theorem \ref{thm_alesker_fourier}). 

\begin{Theorem} \label{thm_main_asymptotic}
As formal power series, 
\begin{align*} 
\sum_{k=0}^\infty \dim \Val_k^{\Sp(\infty)} x^k &
=\frac{x^4-3x^3+6x^2-3x+1}{(1-x)^7(1+x)^3}\\
\sum_{k=0}^\infty \dim \Val_k^{\Sp(\infty)\cdot \U(1)} x^k &
=\frac{x^6-2x^5+2x^4+2x^2-2x+1}{(x^2+1)(x^2+x+1)(1+x)^2(1-x)^6}\\
 \sum_{k=0}^\infty \dim \Val_k^{\Sp(\infty)\cdot \Sp(1)} x^k &
=\frac{x^5+2x^4+x^3+1}{(x^2+1)(x^2+x+1)(1+x)^2(1-x)^4}.
\end{align*}
\end{Theorem}

Explicit formulas and some numerical values for small $k$ can be found at the end of Section \ref{sec_global}. 

In the $\U(n)$-case, \eqref{eq_global_restriction_un},
\eqref{eq_poincare_series_un} and the Alesker-Fourier transform (see
Theorem \ref{thm_alesker_fourier}) are sufficient to compute the
dimension of $\Val_k^{\U(n)}$ for all $k$. In the $\spn$-case, however, we only
get the dimensions in the ranges $0 \leq k \leq n$ and $3n \leq k \leq 4n$. The third main theorem
closes this gap, but the resulting formula is of a combinatorial type and far
from being closed. 

In order to state this theorem, we make the trivial but extremely useful observation that the right action of $\Sp(1)
\cong \SU(2)$ induces the structure of a $\SU(2)$-module on $\Val^{\spn}$. Since
irreducible representations of $\SU(2)$ are indexed by integers,
we have a decomposition 
\begin{equation} \label{eq_decomposition_valk}
 \Val_k^{\spn} \cong \bigoplus_{l=0}^\infty m_k^l V_l,
\end{equation}
where $V_l$ is the unique irreducible $\SU(2)$-representation of dimension $l+1$
and where the coefficients $m_k^l$ are natural numbers. Knowing these
coefficients is equivalent to knowing the character of $\Val_k^{\spn}$, which is
an element in the ring of Laurent polynomials $\mathbb{Z}[t,t^{-1}]$. Also, we
can
read off the dimensions of the spaces $\Val_k^{\spn}$, $\Val_k^{\spnu}$ and
$\Val_k^{\spnsp}$ from the decomposition \eqref{eq_decomposition_valk}: 

\begin{align} \label{eq_dim_sp_abstract}
 \dim \Val_k^{\Sp(n)} & =\sum_l (l+1)m^l_k\\
 \dim \Val_k^{\Sp(n) \cdot \U(1)} & =\sum_{l \equiv 0 (2)} m^l_k \label{eq_dim_spnu1_abstract} \\
\dim \Val_k^{\Sp(n) \cdot \Sp(1)} & =m_k^0. \label{eq_dim_spnsp1_abstract}
\end{align}

In order to describe our third main theorem, we recall some terminology
and conventions and refer to \cite{fulton97} and Section \ref{sec_back_rep} for
details. A {\it Young
diagram} is an arrangement of a finite number of boxes into rows, with less
boxes in lower rows, all aligned to the left. It can be uniquely described by the tuple
$\lambda=(\lambda_1,\ldots,\lambda_k)$ with $\lambda_1 \geq \lambda_2 \geq
\ldots \geq \lambda_k > 0$ where $\lambda_i$ is the number of boxes in the
$i$-th row. Then $k$ is called the {\it depth} of $\lambda$. The {\it weight}
$|\lambda|$ is the number of boxes in $\lambda$. We will call a Young diagram
{\it even}
if the number of boxes in each row is even.

To each Young diagram $\lambda$ there is a corresponding {\it Schur polynomial}
$s_\lambda$ whose definition will be recalled in Section \ref{sec_back_rep}. 

For $n \geq 1$ and $m \geq 0$, let us define the polynomials (over the ring
$\mathbb{Z}[t,t^{-1}]$ of
Laurent polynomials in $t$)
\begin{align}
 E_n(x) & :=\sum_\lambda s_\lambda(tx,t^{-1}x),\label{eq_def_en}\\
F_{n,m}(x) & := \sum_\lambda s_\lambda(tx,t^{-1}x,t,t^{-1}) \label{eq_def_fnm}
\end{align} 
where the first sum is over all even Young diagrams $\lambda$ with
$\lambda_1 \leq 2n$ and the second sum is over all even Young diagrams $\lambda$
with $\lambda_1 \leq 2n$ and $|\lambda|=2m$.

\begin{Theorem} \label{thm_main_thm}
The characters of the $\SU(2)$-representations $\Val_k^{\spn}$ satisfy 
\begin{align*}
\sum_{k=0}^{4n} \ch(\Val_k^{\Sp(n)})x^k & = E_n(x)
-F_{n-1,2n}(x)-(1+x(t^4+t^{-4})+x^2)F_{n-1,2n-1}(x)\\
& \quad +x(1+ x(t^4+1+t^{-4})+ x^2)F_{n-1,2n-2}(x).
\end{align*}
\end{Theorem}

From the characters, one can easily obtain the dimensions of the spaces
of $\spn, \spnu, \spnsp$-invariant valuations (see \eqref{eq_dim_sp_abstract},
\eqref{eq_dim_spnu1_abstract}, \eqref{eq_dim_spnsp1_abstract}). The values for small dimensions $n$ can be computed with the help of Theorem \ref{thm_main_thm} and a computer algebra system, a list is given at the end of Section \ref{sec_dimformula}.

\subsection*{Plan of the paper}
In the next section, we will collect some known facts from the theory of convex
valuations. The only new statement is
Proposition \ref{prop_exact_seq}, which will be central in the proof of Theorems
\ref{thm_main_asymptotic} and \ref{thm_main_thm}. 

In Section \ref{sec_back_rep}, we give some background on quaternionic vector spaces and quaternionic groups and on representation
theory for the groups $\SU(2)$ and $\GL(n,\C)$. 

In Section \ref{sec_invthe_spn} we collect some facts on invariants under the group $\Sp(n)$ which will be used later on. 

The proof of Theorem \ref{thm_local_global} is contained in
Section \ref{sec_global}. It uses several tools from the theory of convex
valuations: normal cycle, Klain's embedding theorem and the Alesker-Fourier
transform. Theorem \ref{thm_main_asymptotic} is proved in the same section. 

In Section \ref{sec_dimformula}, we use some computations for
$\SU(2)$-representations to prove Theorem \ref{thm_main_thm}. Finally, in the
appendix, we will give the rather technical proof of a lemma which is
used in Section \ref{sec_global}. 

\subsection*{Acknowledgments} I wish to thank Semyon Alesker, Joseph Fu and Franz Schuster for useful discussions and comments on this paper. 

\section{Background from valuation theory}
\label{sec_back_val}

We refer to \cite{bernig_aig10} for a recent survey on the subject of {\it algebraic integral geometry}, to which the present paper gives a contribution. 

In this section, $V$ will be a finite dimensional oriented Euclidean vector
space of dimension $m$. The set of compact convex subsets in $V$ is
denoted by $\K(V)$, it is endowed with the Hausdorff metric. 

\begin{Definition}
A (convex) valuation is a map $\mu:\K(V) \to \C$ which is finitely
additive in the following sense: 
 \begin{displaymath}
  \mu(K \cup L)+\mu(K \cap L)=\mu(K)+\mu(L)
 \end{displaymath}
whenever $K,L, K \cup L \in \K(V)$. 
\end{Definition}

The space of continuous, translation invariant valuations is denoted by $\Val$.
With the topology of uniform convergence on compact subsets, it is a Fr\'echet
space, whose dimension is infinite (if $m \geq 2$).  

The main examples of continuous translation invariant valuations are the volume, the Euler characteristic, the intrinsic
volumes and more generally mixed volumes. 

The space $\Val$ comes with a natural grading found by McMullen \cite{mcmullen77}.
An element $\mu \in \Val$ is said to be {\it homogeneous} of degree $k$ if
$\mu(tK)=t^k \mu(K)$ for all $t \geq 0$. $\mu$ is called {\it even} if
$\mu(-K)=\mu(K)$ and {\it odd} if $\mu(-K)=-\mu(K)$ for all $K$. 
The space of $k$-homogeneous even/odd valuations is denoted by $\Val_k^\pm$.

The McMullen grading is given by 
\begin{displaymath}
 \Val=\bigoplus_{\substack{k=0,\ldots,m\\ \epsilon=\pm}} \Val_k^\epsilon. 
\end{displaymath}

For a subgroup $G$ of the rotation group $\SO(V)$, $\Val^G$ denotes the subspace
of $\Val$ of $G$-invariant valuations. If $G$ acts transitively on the unit sphere, then $\Val^G$ is
finite-dimensional \cite{alesker_survey07}. The list of connected closed
subgroups of $\SO(m)$ with this property was given in the introduction (see  \eqref{eq_series} and \eqref{eq_exceptions}). If
$G$ is such a group, then $\Val^G$ consists only of even valuations
\cite{bernig_g2} and we may apply {\it Klain's embedding theorem} \cite{klain00}.
 
More precisely, let $\mu \in
\Val_k^+$. For a $k$-dimensional subspace $E \subset V$, the restriction of
$\mu$ to $\K(E)$ is a multiple of the $k$-dimensional volume in $E$, i.e. 
\begin{displaymath}
 \mu(K)=\kl_\mu(E) \vol_k(K), \quad K \in \K(E)
\end{displaymath}
for some number $\kl_\mu(E)$. The function $\kl_\mu:\Gr_k(V) \to \C$ (where
$\Gr_k$ is the Grassmannian of $k$-dimensional subspaces in $V$) is called the
{\it Klain function} of $\mu$. 

\begin{Theorem} \label{thm_klain}
 The map
\begin{align*}
 \kl:\Val_k^+ & \to C(\Gr_k)\\
\mu & \mapsto \kl_\mu
\end{align*}
is injective. 
\end{Theorem}

Alesker defined the important dense subspace of smooth translation invariant
valuations $\Val^{sm} \subset \Val$. Its definition will be given below. 
\begin{Theorem} \label{thm_alesker_fourier}
 Let $\mu \in \Val_k^{+,sm}$. Then there exists a unique valuation $\F \mu \in
\Val_{m-k}^{+,sm}$ with 
\begin{displaymath}
 \kl_{\F \mu}(E)=\kl_\mu(E^\perp), \quad \forall E \in \Gr_{m-k}. 
\end{displaymath}
\end{Theorem}

The map $\F:\Val^{+,sm} \to \Val^{+,sm}$ is called {\it Alesker-Fourier
transform}. It can be extended to odd valuations \cite{alesker_fourier} but we will not need
this here. Other notations for the Alesker-Fourier transform of $\mu$ are
$\mathbb{D}\mu$ and $\hat \mu$. 

We next describe an important construction of translation invariant valuations,
called {\it normal cycle map} \cite{zaehle86}. In order to simplify notation, we use
the following convention: 

{\it All differential forms are
assumed to be complex-valued.}

By $SV=V \times S(V)$ we denote the unit sphere bundle of $V$. The product
structure on $SV$ induces a bi-grading on the space $\Omega^*(SV)$ of smooth
differential forms, and we denote by $\Omega^{k,l}(SV)$ the space of forms of
bi-degree $(k,l)$. The subspace of translation invariant forms is denoted by a superscript $tr$, i.e. $\Omega^{k,l}(SV)^{tr}=\Lambda^kV^* \otimes \Omega^l(S(V))$. If a group $G$ acts on $V$, then the subspace of $\Omega^{k,l}(SV)$ of translation and $G$-invariant elements is denoted by a superscript $\bar G$. 

On $SV$, there is
a canonical $1$-form $\alpha$ defined by 
\begin{displaymath}
\alpha|_{(x,v)}(w)=\langle v,d\pi(w)\rangle,
\end{displaymath}
where $\pi:SV \to V$ is the canonical projection. The kernel of $\alpha$
defines a contact distribution $Q:=\ker \alpha$, i.e. $SV$ is a $2m-1$-dimensional {\it contact manifold}. 

The {\it Reeb vector field} $T$ is
defined by $T_{(x,v)}=(v,0)$, note that $\alpha(T)=1$. At each point $(x,v)$,
$Q_{(x,v)}$ is the orthogonal sum of two copies of $T_v S(V)$ and we
have an orthogonal splitting 
\begin{displaymath}
T_{(x,v)}SV=\mathbb{R} \cdot T_{(x,v)} \oplus T_v S(V) \oplus T_v S(V).
\end{displaymath}

To any $K \in \K(V)$, we can associate its {\it
normal cycle} $\nc(K)$. The normal cycle is an $m-1$-dimensional Federer-Fleming
current in $SV$. Its support
is the set of pairs $(x,v)$, where $x \in \partial K$ and $v$ is an
outer normal vector to $K$ at $x$. The current $\nc(K)$ is a cycle (i.e.
$\partial \nc(K)=0$) and Legendrian (i.e. $\nc(K) \llcorner \alpha=0$).
 
Let $0 \leq k < m$ and let $\omega \in \Omega^{k,m-1-k}(SV)^{tr}$. We define a continuous translation invariant valuation of degree $k$ by setting 
\begin{displaymath}
 K \mapsto \int_{\nc(K)} \omega.
\end{displaymath}
Linear combinations of valuations of this type and of the volume are
called {\it smooth} (see \cite{alesker_val_man2} for equivalent definitions). The dense subspace of $\Val$ of all smooth valuations is denoted by $\Val^{sm}$. 

An important
fact (see \cite{alesker_survey07}) is that if $G$ is a subgroup of $\SO(V)$ acting transitively on the unit
sphere, then 
\begin{equation} \label{eq_inclusion_valg_smooth}
 \Val^G \subset \Val^{sm}.
\end{equation}

The normal cycle map may be seen as a surjective map 
\begin{align}
\nc: \Omega^{k,m-1-k}(SV)^{tr} & \twoheadrightarrow \Val_k^{sm} \nonumber\\
\omega & \mapsto \left(K \mapsto \int_{\nc(K)} \omega\right).
\label{eq_nc_onto} 
\end{align}

This map is clearly not injective, since vertical forms (i.e. multiples of the
contact form) and exact forms are in the kernel. The kernel was described in
\cite{bernig_broecker07} in terms of the {\it Rumin operator} $D$ from
\cite{rumin94}. This second order differential
operator can be defined on any contact manifold, in particular on the
sphere bundle $SV$. 

Given $\omega \in \Omega^{k,m-1-k}(SV)^{tr}$, there exists a
unique vertical form $\xi$ such that $d(\omega+\xi)$ is vertical and
$D\omega:=d(\omega+\xi)$. 

It follows immediately that $D$ vanishes on closed forms $\omega$ (since
$d\omega=0$ is vertical), on vertical forms (take $\xi:=-\omega$) and on
multiples of $d\alpha$: 
\begin{displaymath}
 D (d\alpha \wedge \tau)=d(d\alpha \wedge \tau + \alpha \wedge d\tau)=0.
\end{displaymath}

The next theorem is a special case of the main theorem in
\cite{bernig_broecker07}. 

\begin{Theorem} \label{thm_bernig_broecker}
 $\omega \in \ker \nc$ if and only if $D\omega=0$ and
$\pi_*\omega=0$. 
\end{Theorem}

Note that the second condition is always satisfied if $k>0$. 

Using arguments as in Section 2 in Rumin's paper 
\cite{rumin94}, we show that $\Val_k^{sm}$ fits into some exact sequence.

For this, we define the following spaces. 
\begin{align*}
\mathcal{I}^{k,l}(SV)^{tr}  & :=\{ \omega \in \Omega^{k,l}(SV)^{tr}: \omega=
\alpha \wedge
\xi+ d\alpha \wedge \psi, \\
& \quad \xi \in \Omega^{k-1,l}(SV)^{tr}, \psi \in
\Omega^{k-1,l-1}(SV)^{tr}\}\\
\Omega_v^{k,l}(SV)^{tr} & :=\{\omega \in \Omega^{k,l}(SV)^{tr}: \alpha \wedge
\omega=0\}\\
\Omega_h^{k,l}(SV)^{tr} & :=\Omega^{k,l}(SV)^{tr} / \Omega_v^{k,l}(SV)^{tr}\\
\Omega_p^{k,l}(SV)^{tr} & :=\Omega^{k,l}(SV)^{tr} / \mathcal{I}^{k,l}(SV)^{tr}.
\end{align*}

Multiplication by the symplectic form $-d\alpha$ induces an operator
$L:\Omega_h^{k,l}(SV)^{tr} \to \Omega_h^{k+1,l+1}(SV)^{tr}$ which is an
injection for
$k+l \leq m-2$, moreover
$\Omega_p^{k,l}(SV)^{tr} \simeq \Omega_h^{k,l}(SV)^{tr}/L\Omega_h^{k-1,l-1}(SV)^{tr}$.
The
exterior derivative induces an operator $d_Q:\Omega_p^{k,l}(SV)^{tr} \to
\Omega_p^{k,l+1}(SV)^{tr}$. 

The Rumin operator vanishes on $\mathcal{I}^{k,m-k-1}(SV)^{tr}$, hence it
induces an operator (which we denote by the same letter)
$D:\Omega_p^{k,m-k-1}(SV)^{tr} \to \Omega_v^{k,m-k}(SV)^{tr}$. 

\begin{Lemma} \label{lemma_translation_poincare}
Let $\omega \in \Omega^{k,l}(SV)^{tr}$ with $d\omega=0$. 
\begin{enumerate}
\item In the case $0<l<m-1$, there exists $\phi \in \Omega^{k,l-1}(SV)^{tr}$
with
$d\phi=\omega$. 
\item In the case $l=0$, $\omega \in \Lambda^kV^* \otimes \C$.
\item In the case $k=0$, $l=m-1$, there exists $\phi \in
\Omega^{0,m-2}(SV)^{tr}$ with $d\phi=\omega$ provided that $\pi_* \omega=0$. 
\end{enumerate}
\end{Lemma}

\proof
We write 
\begin{displaymath}
\omega=\sum_{i=1}^q c_i \kappa_i \wedge \tau_i
\end{displaymath}
where $\kappa_1,\ldots,\kappa_q$ is a basis of $\Lambda^kV^*$ and
$\tau_1,\ldots,\tau_q$ are $l$-forms on the unit sphere $S(V)$. Then 
\begin{displaymath}
d\omega=(-1)^k \sum_{i=1}^q c_i \kappa_i \wedge d\tau_i
\end{displaymath}
which shows that all $\tau_i$ are closed. 

If $0 < l < m-1$, then $H_{dR}^l(S^{m-1})=0$ and we find $\rho_i \in
\Omega^{l-1}(S(V))$ with $d\rho_i=\tau_i$. Then $\phi:=(-1)^k \sum_{i=1}^q c_i
\kappa_i \wedge \rho_i$ satisfies $d\phi=\omega$. If $l=0$, then all $\tau_i$
are constant and hence $\omega \in \Lambda^kV^* \otimes \C$.

The last statement follows from the fact that
$\Omega^{0,m-1}(SV)^{tr}=\Omega^{m-1}(S(V))$ and that $H^{m-1}_{dR}(S(V))$ is
one-dimensional.  
\endproof

\begin{Proposition} \label{prop_exact_seq}
Let $0 \leq k \leq m$. There is an exact sequence 
\begin{displaymath}
0 \to \Lambda^k V^* \otimes \C \hookrightarrow \Omega^{k,0}_p(SV)^{tr}
\stackrel{d_Q}{\to}
\Omega^{k,1}_p(SV)^{tr} \stackrel{d_Q}{\to} \ldots \stackrel{d_Q}{\to}
\Omega^{k,m-k-1}_p(SV)^{tr} \stackrel{\nc}{\to} \Val_k \to 0
\end{displaymath}
\end{Proposition}

The proof follows the arguments in \cite{rumin94}, with the Poincare lemma
replaced by Lemma \ref{lemma_translation_poincare} above.  

\proof
 Since $\Val_m$ is spanned by the Lebesgue measure, there is an exact
sequence 
\begin{displaymath}
 0 \to \Lambda^m V^* \otimes \C \to \Val_m \to 0,
\end{displaymath}
which is the case $k=m$ in the statement. 

Let us suppose that $k<m$. It is clear that the sequence is closed. Let us check that it is exact.

\begin{itemize}
\item  Since
$\Omega^{k,0}_p(SV)^{tr}=\Omega^{k,0}(SV)^{tr}$, the injectivity on
the left hand side is trivial. 
\item For $0 \leq l<m-k-1$, let $\omega \in \Omega^{k,l}(SV)^{tr}$ with
$d\omega=\alpha
\wedge \xi + d\alpha \wedge \psi \in \mathcal{I}_{k,l+1}(SV)^{tr}$. 

Letting $\omega':=\omega-\alpha \wedge \psi \in \Omega^{k,l}(SV)^{tr}$ and
$\xi':=\xi+d\psi$, we obtain $d \omega'=\alpha \wedge \xi'$. Differentiating
yields $0=d\alpha \wedge \xi' - \alpha \wedge d\xi'$. In other words,
$L(\xi'|_Q)=0$. By the injectivity of $L$ in degree $k+l<m-1$, we get
$\xi'|_Q=0$ which implies that $d\omega'=0$. 

If $l > 0$, Lemma \ref{lemma_translation_poincare} implies that there exists
$\phi \in \Omega_{k,l-1}(SV)^{tr}$ with $d\phi=\omega'$. Hence
$[\omega]=[\omega']=[d\phi]=d_Q[\phi]$, which shows that $[\omega]$ is
$d_Q$-exact.  
\item In the same situation, if $l=0$, then $\omega'$ is a translation invariant $k$-form on
$V$, hence $[\omega]=[\omega']$ is in the image of $\Lambda^k V^* \otimes \C$.  
\item The map $\nc$ on the right hand side is surjective by \eqref{eq_nc_onto}. 
\item If $[\omega] \in \Omega_p^{k,m-k-1}(SV)^{tr}$ lies in the kernel of $\nc$,
then
$D[\omega]=d(\omega+\xi)=0$ for some form $\xi \in
\Omega_v^{k,m-k-1}(SV)^{tr}$. Then 
$\omega':=\omega+\xi$ is a closed translation invariant form of
bi-degree $(k,m-k-1)$. If $k>0$, then, by Lemma \ref{lemma_translation_poincare},
there exists
$\phi \in \Omega^{k,m-k-2}(SV)^{tr}$ with $d\phi=\omega'$. Hence
$[\omega]=[\omega']=[d\phi]=d_Q[\phi]$ is $d_Q$-exact. 
\item In the same situation, if $k=0$, then
$\pi_* \omega'=\pi_*\omega=0$ by Theorem \ref{thm_bernig_broecker} (note that
$\pi_*\xi=0$, since $\xi$ is vertical). We may thus apply Lemma
\ref{lemma_translation_poincare} to find $\phi \in \Omega^{0,m-2}(SV)^{tr}$ with
$d\phi=\omega'$. Hence $[\omega]$ is $d_Q$-exact. 
\end{itemize}
\endproof 

\begin{Corollary} \label{cor_exact_seq_group}
Let $G$ be a closed subgroup of $\SO(V)$ acting transitively on the unit sphere
of $V$. Let $0 \leq k \leq m$. There is an exact sequence 
\begin{displaymath}
0 \to (\Lambda^k V^* \otimes \C)^G \hookrightarrow \Omega^{k,0}_p(SV)^{\bar G}
\stackrel{d_Q}{\to}
\Omega^{k,1}_p(SV)^{\bar G} \stackrel{d_Q}{\to} \ldots \stackrel{d_Q}{\to}
\Omega^{k,m-k-1}_p(SV)^{\bar G} \stackrel{\nc}{\to} \Val_k^G \to 0
\end{displaymath}
\end{Corollary}

\proof
In Lemma \ref{lemma_translation_poincare}, if $\omega$ is
$G$-invariant, then $\phi$ may be chosen $G$-invariant too (just average
$\phi$ with respect to the Haar measure on $G$). The rest of the proof is
analogous. 
\endproof
\section{Background from combinatorics and representation theory}
\label{sec_back_rep}

\subsection{Young diagrams and Schur functions}
\label{subsec_young_diagrams}

Our main reference for Young diagrams is \cite{fulton97}.

Let $\lambda=(\lambda_1,\ldots,\lambda_k)$ be a given partition, i.e. $\lambda_1
\geq \lambda_2 \geq \ldots \lambda_k \geq 0$. The {\it Young diagram}
associated to $\lambda$ has $\lambda_i$ boxes in the $i$-th row, all aligned
to the left. For example, if $\lambda=(4,1,1)$, the Young diagram is given by  

\begin{displaymath}
\lambda=\yng(4,1,1)
\end{displaymath}

The {\it depth} of a Young diagram is the number of its rows, the {\it weight}
is the number of its boxes. The {\it conjugate} $\tilde \lambda$ of a Young diagram $\lambda$ is obtained
by reflecting it at the diagonal. In our example
\begin{displaymath}
\tilde \lambda=\yng(3,1,1,1)
\end{displaymath}

The Young diagram $\lambda$ {\it dominates} another Young diagram $\mu$ if
\begin{displaymath}
 \sum_{i=1}^m \lambda_i \geq \sum_{i=1}^m \mu_i
\end{displaymath}
for all $m$. In this case, we write $\lambda \trianglerighteq \mu$.

A {\it semi-standard tableaux} on $\lambda$ is given by putting one of the
numbers $1,2,\ldots,m$ (where $m$ may be different from the weight of
$\lambda$) into each box of $\lambda$ in such a way that the numbers weakly
increase in each row from left to right and strictly increase in each column
from top to bottom. 

For instance the semi-standard
tableaux on $\lambda=\yng(2,1)$ with $m=3$ are 
\begin{multline*}
 \young(11,2) \quad \young(12,2) \quad \young(13,2) \quad
\young(11,3) \quad \young(12,3) \quad \young(13,3) \quad \young(22,3) \quad
\young(23,3)
\end{multline*}

Let $\lambda$ be a Young diagram and let
$\mu=(\mu_1,\ldots,\mu_m)$ be a multi-index. The number
of semi-standard tableaux on $\lambda$ with exactly $\mu_1$ 1's, $\mu_2$ 2's and
so on is denoted by $K_{\lambda \mu}$ and is called {\it Kostka number}. A
non-trivial fact is that $K_{\lambda \mu}$ does not depend on the order of
$\mu$, hence we may assume that $\mu$ is itself a Young diagram. Since in any
semi-standard tableaux, the numbers in the boxes of the $i$-th row are at
least $i$, we have
$K_{\lambda \mu}=0$ unless $\lambda \trianglerighteq \mu$. 

For a given Young diagram $\lambda$ and a number $m$, the
polynomial 
\begin{displaymath}
s_\lambda(x_1,x_2,\ldots,x_m):=\sum_\mu K_{\lambda \mu} x_1^{\mu_1} x_2^{\mu_2}
\ldots x_m^{\mu_m},
\end{displaymath}
where $\mu$ ranges over all multi-indices, is called
{\it Schur
polynomial} of $\lambda$. By convention, if $\lambda=(0,\ldots,0)$ is the
empty Young diagram, then $s_\lambda(x_1,\ldots,x_m)=1$.  

In our example, $\lambda=(2,1)$ and $m=3$, and we get 
\begin{displaymath}
 s_\lambda(x_1,x_2,
x_3)=x_1^2x_2+x_1x_2^2+2x_1x_2x_3+x_1^2x_3+x_1x_3^2+x_2^2x_3+x_2x_3^2.
\end{displaymath}
In general,
$s_\lambda(x_1,\ldots,x_m)$ is a symmetric polynomial. Moreover, the set of
all Schur polynomials $s_\lambda(x_1,\ldots,x_m)$, as $\lambda$ ranges over all
Young diagrams of depth $\leq m$, is a basis of the vector space of all symmetric polynomials in
$x_1,\ldots,x_m$. 

\subsection{Representations of $\SU(2)$}

Recall that a maximal torus in $\SU(2)$ is $\U(1)=S^1$, which we embed by $z
\mapsto \left(\begin{array}{c c} z & 0\\0 & z^{-1}\end{array}\right)$ (see
also \eqref{eq_embedding_h}). If $(V,\phi)$ is a finite-dimensional
representation of $\SU(2)$, then 
\begin{displaymath}
 V=\bigoplus_\alpha V_\alpha,
\end{displaymath}
where $\alpha$ ranges over all weights of $V$ (which are integers) and 
\begin{displaymath}
 V_\alpha=\{v \in V| \phi(z)(v)=z^\alpha v\}
\end{displaymath}
is the corresponding {\it weight space}. A non-zero vector $v \in V_\alpha$ is called
{\it vector of weight $\alpha$}.

The {\it character} of $(V,\phi)$ is defined by 
\begin{displaymath}
 \ch(V)=\sum_{\alpha=-\infty}^\infty \dim V_\alpha t^\alpha \in
\mathbb{Z}[t,t^{-1}].
\end{displaymath}
 
Let $V \cong \C^2$ be the standard representation of $\SU(2)$ and set
$V_k:=\Sym^k V$. Then $V_k$ is an irreducible
$\SU(2)$-representation of dimension $k+1$. Any finite-dimensional irreducible
representation of $\SU(2)$ is isomorphic to $V_k$ for some $k$. We have 
\begin{displaymath}
 \ch(V_k)=\sum_{j=0}^k t^{2j-k} \in \mathbb{Z}[t,t^{-1}].
\end{displaymath}

We will frequently use the following rules, whose proofs can be found in
\cite{fulton_harris91}:
\begin{align}
V_k \otimes V_l & \simeq V_{k+l} \oplus V_{k+l-2} \oplus \ldots \oplus V_{|k-l|}
\quad \text{(Clebsch-Gordan rule)} \label{eq_clebsch_gordan} \\
 \Sym^k V_2 & \simeq \bigoplus_{l=0}^{\lfloor
k/2\rfloor} V_{2k-4l}. \label{eq_symmetric_powers}
\end{align}

It will be convenient to state many formulas in the representation ring $R\SU(2)$, which is the
free $\mathbb{Z}$-module generated by variables $V_0,V_1,V_2,\ldots$, with
multiplication given by the Clebsch-Gordan rule. 

The spaces of differential forms which are considered in this paper
are {\it graded commutative $\SU(2)$-algebras}, i.e.
graded $\SU(2)$-modules $W$ endowed with a graded-commutative product
map $W \otimes W \to W$ which is a $\SU(2)$-morphism and which is compatible
with the grading.

\subsection{Representations of $\GL(n,\C)$}

Irreducible complex representations of $\GL(n,\C)$ are classified
by sequences $\lambda=(\lambda_1,\ldots,\lambda_n)$ of integers with $\lambda_1 \geq \lambda_2 \geq \ldots \geq \lambda_n$. The standard ($n$-dimensional) representation corresponds to $\lambda=(1)$, while the $1$-dimensional determinantal representation corresponds to $\lambda=(1,1,\ldots,1)$.  

An explicit description of the
representation $\Gamma_\lambda$ corresponding to an arbitrary sequence $\lambda$ can be
found in any standard text on representation theory, e.g.
\cite{fulton_harris91} or \cite{fulton97}. If $\lambda_n \geq 0$, then the character of $\Gamma_\lambda$ equals the Schur function $s_\lambda(x_1,\ldots,x_n)$. 

Let $V$ be the standard representation of $\GL(n,\C)$. A standard fact from
representation theory of the general linear group (\cite{fulton97}, \S 8.3)
tells us that 
\begin{equation} \label{eq_dec_prod_alt}
 \Lambda^{k_1}V \otimes \ldots \otimes \Lambda^{k_s}V \simeq \bigoplus_\nu K_{\tilde \nu \mu} \Gamma_\nu,
\end{equation}
where $\mu=(k_1,\ldots,k_s)$ and $\nu$ ranges over all Young diagrams of
depth $\leq n$, $\tilde \nu$ is the
conjugate of $\nu$ and $K_{\tilde \nu \mu}$ is the Kostka 
number from Subsection \ref{subsec_young_diagrams}.   

\subsection{Symplectic groups}
\label{subsec_symp_groups}

Let $\h$ be the $4$-dimensional real vector space of quaternions. This space has
an
$\R$-basis given by $1,i,j,k$. The algebra structure of $\h$ is defined by
$i^2=j^2=k^2=-1, ij=-ji=k, jk=-kj=i, ki=-ik=j$. The conjugate of a quaternion
$q=a+bi+cj+dk$ is the quaternion $\bar q:=a-bi-cj-dk$. The norm of
$q$ is the real number $q \bar q$. 

We will consider $\h$ as a complex vector space, i.e. we write 
\begin{equation} \label{eq_dec_h}
 \h=\C \oplus j\C. 
\end{equation}

Multiplication from the left on $\h=\C \oplus j\C=\C^2$ defines an embedding of
$\h$ into $M_2\C$. Explicitly, 
\begin{align}
\h & \hookrightarrow M_2\C \nonumber \\
 z_1 + jz_2 & \mapsto \left(\begin{array}{c c}z_1 & -\bar z_2\\z_2 & \bar
z_1\end{array}\right) \label{eq_embedding_h} 
\end{align}

The quaternions of norm $1$ form a subgroup $\Sp(1)$ of $\h^*$. It contains the subgroup $\U(1)$ of all complex numbers of norm $1$ (i.e. quaternions of the form $a+ib$ with $a^2+b^2=1$). The embedding
\eqref{eq_embedding_h} identifies $\Sp(1)$ with $\SU(2)$, which is homeomorphic
to a three-dimensional sphere.

Let $V$ be a quaternionic (right) vector space of dimension $n$. We endow $V$
with a quaternionic Hermitian form $K$, i.e. an $\R$-bilinear form 
\begin{displaymath}
 K:V \times V \to \h
\end{displaymath}
such that  
\begin{enumerate}
\item $K$ is conjugate $\h$-linear in the first and $\h$-linear in the second
factor, i.e. 
\begin{displaymath}
 K(v q, w r)=\bar q K(v,w) r \quad q,r \in \h,
\end{displaymath}
\item $K$ is Hermitian in the sense that 
\begin{displaymath}
 K(w,v)=\overline{K(v,w)},
\end{displaymath}
\item $K$ is positive definite, i.e. 
\begin{displaymath}
 K(v,v)>0 \quad \forall v \neq 0. 
\end{displaymath}
\end{enumerate}

The standard example of such a form is given by
\begin{displaymath}
V=\h^n, K(v,w)=\sum_{i=1}^n \bar v_i w_i, \quad v=(v_1,\ldots,v_n),
w=(w_1,\ldots,w_n) \in \h^n.
\end{displaymath}
 
Recall that $\GL(V,\h)=\GL(n,\h)$ is the group of $\h$-linear automorphisms of $V$. The subgroup of $\GL(V,\h)$ of all elements preserving $K$ is called {\it compact
symplectic group} and denoted by $\Sp(V,K)$ or $\Sp(n)$. It acts from the left on
$V$. 

With respect to the decomposition \eqref{eq_dec_h}, we can decompose 
\begin{displaymath}
 K(v,w)=H(v,w)+jQ(v,w).
\end{displaymath}
Then $H$ is a complex Hermitian form and $Q$ a skew-symmetric complex linear
form. Moreover, $Q(v,w)=H(vj,w)$ and $H(vj,wj)=\overline{H(v,w)}$. If $\U(2n)$
denotes the unitary group with respect to $H$, and $\Sp_{2n}\C$ the symplectic
group with respect to $Q$ (i.e. the subgroup of $\GL(2n,\C)$ consisting of all
elements preserving $Q$) then 
\begin{displaymath}
 \Sp(n)=\U(2n) \cap \Sp_{2n}\C.
\end{displaymath}
Hence $\Sp(n)$ is a compact form of $\Sp_{2n}\C$, a fact which will be used in Section \ref{sec_invthe_spn}.

We let $\Sp(1)$ and its subgroup $\U(1)$ act on $V$ from the right. The actions
of $\Sp(n) \times \Sp(1)$ and $\Sp(n) \times \U(1)$ on $V$ are not
faithful: in both cases the kernel consists of the two elements
$(Id,1)$ and $(-Id,-1)$ and is isomorphic to $\mathbb{Z}_2$. Hence the factor
groups
\begin{align*}
 \Sp(n)\cdot \Sp(1) & := \Sp(n) \times \Sp(1)/ \mathbb{Z}_2\\
 \Sp(n)\cdot \U(1) & := \Sp(n) \times \U(1)/ \mathbb{Z}_2
\end{align*}
act faithfully on $V$. Since $\Sp(n)$ acts transitively on the unit sphere,
the same holds true for $\Sp(n) \cdot \Sp(1)$ and $\Sp(n) \cdot \U(1)$. 

Let $v_0 \in V$ be a unit vector. Let us assume that $v_0$ is the first standard
coordinate vector. Then a pair $(A,z) \in \Sp(n) \times \Sp(1)$ stabilizes $v_0$
if and only if $A$ is of the form 
\begin{displaymath}
 A=\left(\begin{array}{c c} z^{-1} & 0 \\ 0 & A'\end{array}\right), \quad A' \in \Sp(n-1).
\end{displaymath}
It follows that the stabilizer of $\Sp(n) \cdot \Sp(1)$ at $v_0$ is isomorphic to
$\Sp(n-1) \cdot \Sp(1)$.  

The tangent space $T_{v_0}S^{4n-1}$ splits as $U \oplus
\tilde V$, where $\tilde V$ is the quaternionic orthogonal complement of $v_0$
and $U$ is
of dimension $3$. Then the action of $\Sp(n-1)\cdot \Sp(1)$ on $\tilde V \cong
\h^{n-1}$ is the
usual
one (i.e. $\Sp(n-1)$ acts from the left, $\Sp(1)$ from the right), while the
action on $U$ is the adjoint action of $\Sp(1)$ (i.e. $(A',z)$ acts by
multiplication by $z^{-1}$ from the left followed by right multiplication by
$z$).  

\subsection{Spherical representations of $\SL(2n,\C)$}

We will need the following proposition which is a consequence of known facts from representation theory. We refer to \cite{goodman_wallach09}, Section 12,  for the terminology in the proof. 

\begin{Proposition}
Let $\lambda$ be a Young diagram and $\Gamma_\lambda$ the corresponding irreducible representation of $\GL(2n,\C)$. Let $\Gamma_\lambda^{\Sp_{2n}\mathbb{C}}$ be the subspace of $\Sp_{2n}\mathbb{C}$-fixed elements. Then 
 \label{prop_spherical}
\begin{displaymath}
\dim \Gamma_\lambda^{\Sp_{2n}\mathbb{C}}=\left\{\begin{array}{c c} 1 & \text{ if } \tilde \lambda \text{ even}\\ 0 & \text{ otherwise.}\end{array} \right.
\end{displaymath}
\end{Proposition}

\proof
Since $\Sp_{2n}\mathbb{C} \subset \SL(2n,\C)$, the determinantal representation of $\GL(2n,\C)$ restricts to the trivial representation of $\Sp_{2n}\mathbb{C}$. We may thus suppose that $\lambda_{2n}=0$ and consider $\Gamma_\lambda$ as a representation of $\SL(2n,\C)$. 

Since the pair $G=\SL(2n,\C), K=\Sp_{2n}\mathbb{C}$ is a {\it spherical pair}, the dimension of $\Gamma_\lambda^K$ is at most $1$. If $\Gamma_\lambda$ contains a non-zero $K$-fixed vector, $\Gamma_\lambda$ is called {\it spherical}.

By Theorem 12.3.13. of \cite{goodman_wallach09} and the explicit computation in (\cite{goodman_wallach09}, Section 12.3.2, Type AII), we may conclude that $\Gamma_\lambda$ is spherical if and only if $\tilde \lambda$ is even.
\endproof


\section{Invariant theory of $\Sp(n)$}
\label{sec_invthe_spn}

Let us fix a sequence of quaternionic isometric embeddings 
\begin{displaymath}
\mathbb{H}^1 \stackrel{\iota_{12}}{\hookrightarrow}
\mathbb{H}^2 \stackrel{\iota_{23}}{\hookrightarrow} \mathbb{H}^3 \stackrel{\iota_{34}}{\hookrightarrow} \ldots 
\end{displaymath}
Then we get a sequence of restriction maps
\begin{displaymath}
(\Lambda^* \h^* \otimes \C)^{\Sp(1)} \stackrel{\iota_{12}^*}{\longleftarrow}
(\Lambda^*
(\h^2)^* \otimes \C)^{\Sp(2)} \stackrel{\iota_{23}^*}{\longleftarrow} (\Lambda^*
(\h^3)^* \otimes \C)^{\Sp(3)} \stackrel{\iota_{34}^*}{\longleftarrow} \ldots  
\end{displaymath}
and define
\begin{displaymath}
 (\Lambda^* (\h^\infty)^* \otimes \C)^{\Sp(\infty)}:=\lim_{\longleftarrow}
(\Lambda^*
(\h^n)^* \otimes \C)^{\Sp(n)}.
\end{displaymath}

The grading, algebra structure and $\SU(2)$-action on each of the spaces
$(\Lambda^* (\h^n)^* \otimes \C)^{\Sp(n)}$ induce the structure of a graded
$\SU(2)$-algebra on $(\Lambda^*(\h^\infty)^* \otimes \C)^{\Sp(\infty)}$.

\begin{Proposition} \label{prop_free_algebra_down}
 The graded $\SU(2)$-algebra $(\Lambda^*(\h^\infty)^* \otimes \C)^{\Sp(\infty)}$ is
freely
generated 
by one copy of $V_2$ in degree $2$. The restriction map
\begin{displaymath}
  u_n:(\Lambda^k (\h^\infty)^* \otimes \C)^{\Sp(\infty)} \to  (\Lambda^k
(\h^n)^* \otimes \C)^{\Sp(n)}
\end{displaymath}
commutes with the $\SU(2)$-action, is surjective for all $n$ and injective for $n \geq k$. 
\end{Proposition}

\proof
Let $V$ be an $n$-dimensional quaternionic vector space with a quaternionic
Hermitian form $K$ as in the previous subsection. We have
an isomorphism between $\Sp_{2n}\C$-modules
\begin{align} 
 \phi:V & \cong V^*\label{eq_iso_v_v*}\\
v & \mapsto [w \mapsto Q(v,w)]. \nonumber
\end{align}

In the following, we consider $V$ as a complex vector space. To each
(complex) polynomial $f: \underbrace{V \times \ldots \times V}_{m} \to \C$ in the real and
imaginary coordinates, we associate a polynomial $\tilde f:
\underbrace{V \times \ldots \times V}_{m} \times \underbrace{V^* \times \ldots \times V^*}_{m} \to \C$ in
the following way:
\begin{enumerate}
 \item If $f(v_1,\ldots,v_m)=v_{ij}$, then $\tilde
f(v_1,\ldots,v_m,\xi_1,\ldots,\xi_m):=v_{ij}$;
\item If $f(v_1,\ldots,v_m)=\bar v_{ij}$, then
$\tilde f(v_1,\ldots,v_m,\xi_1,\ldots,\xi_m):=\xi_{ij}$;
\item $f  \mapsto \tilde f$ is an algebra homomorphism. 
\end{enumerate}

We get another polynomial $\hat{f}:\underbrace{V \times \ldots \times V}_{2m} \to
\C$ by composing with the isomorphism \eqref{eq_iso_v_v*}., i.e. 
\begin{displaymath}
 \hat{f}(v_1,\ldots,v_m,w_1,\ldots,w_m):=\tilde
f(v_1,\ldots,v_m,\phi(w_1),\ldots,\phi(w_m)).
\end{displaymath}

It follows as in (\cite{spivak5_79}, Chapter 13, Addendum 1), that if $f$ is invariant under $\spn$, then
$\tilde f$ and $\hat{f}$ are invariant under
$\Sp_{2n}\C$. 

Now the first fundamental theorem (FFT) for $\Sp_{2n}\C$
(\cite{fulton_harris91}, Prop. F.13), tells us that $\hat{f}$ is a
polynomial in the {\it basic invariants} 
\begin{displaymath}
 (v_1,\ldots,v_{2m}) \mapsto Q(v_i,v_j) \quad 1 \leq i<j \leq 2m.
\end{displaymath}

Unravelling the identifications, we get that $f$ is a polynomial in the three basic invariants 
\begin{displaymath}
 (v_1,\ldots,v_m) \mapsto \left\{ \begin{array}{c c} Q(v_i,v_j) & 1 \leq
i<j\leq m\\ \overline{Q(v_i,v_j)} & 1 \leq i < j \leq m\\
H(v_i,v_j) & 1 \leq i,j \leq m.\end{array}\right.
\end{displaymath}

Since $Q$ and $\bar Q$ are antisymmetric and since $H(v,w)-H(w,v)=2i \im
H(v,w)$, we deduce that the algebra $(\Lambda^*V^*
\otimes \C)^{\spn}$ is
generated by the three basic $2$-forms $Q$, $\bar Q$ and $\im H$. They span
an irreducible $3$-dimensional $\SU(2)$-representation. 

The second fundamental theorem (SFT) for $\Sp_{2n}\C$ (\cite{goodman_wallach09}
Thm. 12.2.15) yields that there are no relations between the basic invariants in
degree less than or equal to $n$, hence the above restriction map is injective
if $k \leq n$. 
\endproof

Let 
\begin{displaymath}
S\mathbb{H}^1 \stackrel{\tilde \iota_{12}}{\hookrightarrow}
S\mathbb{H}^2 \stackrel{\tilde \iota_{23}}{\hookrightarrow} S\mathbb{H}^3 \stackrel{\tilde \iota_{34}}{\hookrightarrow} \ldots 
\end{displaymath}
be the induced sequence of embeddings. We get a sequence of
restriction maps
\begin{displaymath}
(\Omega^{*,*}_h(S\h^1))^{\Sp(1)} \stackrel{\tilde \iota_{12}^*}{\longleftarrow}
(\Omega^{*,*}_h(S\h^2))^{\Sp(2)} \stackrel{\tilde \iota_{23}^*}{\longleftarrow}
(\Omega^{*,*}_h(S\h^3))^{\Sp(3)} \stackrel{\tilde \iota_{34}^*}{\longleftarrow}
\ldots  
\end{displaymath}
and define 
\begin{displaymath}
(\Omega^{*,*}_h(S\h^\infty))^{\Sp(\infty)}:=\lim_{\longleftarrow}
(\Omega^{*,*}_h(S\h^n))^{\Sp(n)},
\end{displaymath}
which is a bi-graded $\SU(2)$-algebra. 

\begin{Proposition} \label{prop_free_algebra_up}
 The bi-graded $\SU(2)$-algebra $(\Omega^{*,*}_h(S\h^\infty))^{\Sp(\infty)}$ is
generated by five copies of $V_2$ in bi-degrees $(1,0),(0,1),(2,0),(1,1),(0,2)$
and one copy of $V_0$ in bi-degree $(1,1)$. The restriction maps 
\begin{displaymath}
 \tilde u_n:(\Omega^{k,l}_h(S\h^\infty))^{\Sp(\infty)} \to
(\Omega^{k,l}_h(S\h^n))^{\Sp(n)}
\end{displaymath}
are surjective for all $k,l,n$ and injective for $k+l \leq n-1$. 
\end{Proposition}

\proof
The proof is based on the same argument as the proof of Proposition \ref{prop_free_algebra_down}.

Fix a unit vector vector $v_0 \in V$. The orthogonal complement $v_0^\perp
\subset V$ splits as 
\begin{displaymath}
 v_0^\perp=U \oplus \tilde V,
\end{displaymath}
where $\tilde V$ is the quaternionic orthogonal complement of $v_0$ and $U:=v_0
\cdot
\h \cap v_0^\perp$ is three-dimensional. The stabilizer of $\Sp(n) \cong \Sp(V,K)$
at $v_0$ is $\Sp(n-1) \cong \Sp(\tilde V,K|_{\tilde V})$, it acts trivially on
$U$. As
$\SU(2)$-representation $U \simeq V_2$, the adjoint representation (compare the
discussion in Subsection \ref{subsec_symp_groups}).

The horizontal hyperplane at the point $(0,v_0) \in SV$ is then given by 
\begin{displaymath}
 Q_{(0,v_0)}=U \oplus \tilde V \oplus U \oplus \tilde V,
\end{displaymath}
and $\Sp(n-1)$ acts trivially on the $U$-factors and diagonally on the second and fourth factor.

As $\SU(2)$-representations we thus have the following isomorphisms.
\begin{align*}
\Omega_h^{*,*}(SV)^{\spn} & \simeq \Lambda^{*,*} Q^*_{(0,v_0)} \otimes \C\\
& \simeq \Lambda^{*,*} (U^* \oplus \tilde V^* \oplus U^* \oplus \tilde V^*)^{\Sp(n-1)}
\otimes \C \\
& \simeq \left(\Lambda^*(U^* \oplus \tilde V^*) \otimes \Lambda^*(U^* \oplus
\tilde V^*)\right)^{\Sp(n-1)} \otimes \C \\
& \simeq \sum_{k,l=0}^3 \Lambda^k U^* \otimes \Lambda^l U^* \otimes
\left(\Lambda^{*-k}\tilde V^* \otimes  \Lambda^{*-l}\tilde V^*\right)^{\Sp(n-1)}
\otimes \C
\end{align*}
Now the algebra $\Lambda^* U^* \otimes \C$ is generated by $U^* \simeq V_2$.
Using arguments as in the previous proof, we get that $\left(\Lambda^{*}\tilde
V^*
\otimes \Lambda^{*}\tilde V^* \otimes \C\right)^{\Sp(n-1)}$ is generated by three
copies of $V_2$ in degrees $(2,0),(1,1)$ and $(0,2)$ and one copy of $V_0$
corresponding to the symplectic form on $SV$. 

These invariant forms are defined on $S\h^n$ for each $n$, behave well
under the restriction maps and can therefore be considered as elements of
$\Omega_h^{*,*}(S\h^\infty)^{\Sp(\infty)}$. The  surjectivity of $\tilde u_n$
follows. 

Since by the SFT for $\Sp_{2n-2}\C$ there are no relations of degree $ \leq n-1$
between the basic $\Sp(n-1)$-invariant forms, the injectivity of $\tilde u_n$ in the case $k+l \leq n-1$
follows.
\endproof

\section{Local and global valuations}
\label{sec_global}

Let us fix again a sequence of quaternionic embeddings 
\begin{displaymath}
\mathbb{H}^1 \stackrel{\iota_{12}}{\hookrightarrow}
\mathbb{H}^2 \stackrel{\iota_{23}}{\hookrightarrow} \mathbb{H}^3 \stackrel{\iota_{34}}{\hookrightarrow} \ldots 
\end{displaymath}

Correspondingly, we have embeddings of the sphere bundles 
\begin{displaymath}
S\mathbb{H}^1 \stackrel{\tilde \iota_{12}}{\hookrightarrow}
S\mathbb{H}^2 \stackrel{\tilde \iota_{23}}{\hookrightarrow} S\mathbb{H}^3
\stackrel{\tilde \iota_{34}}{\hookrightarrow} \ldots 
\end{displaymath}
and restrictions  
\begin{displaymath}
\Val^{\Sp(1)} \stackrel{r_{21}}{\longleftarrow} \Val^{\Sp(2)}
\stackrel{r_{32}}{\longleftarrow} \Val^{\Sp(3)}
\stackrel{r_{43}}{\longleftarrow}
\ldots
\end{displaymath} 

Let 
\begin{displaymath}
 \Val^{\Sp(\infty)}:=\lim_{\longleftarrow} \Val^{\Sp(n)}
\end{displaymath}
be the inverse limit. An element of $\Val^{\Sp(\infty)}$ is thus a
sequence $(\mu_1,\mu_2,\ldots)$ with $\mu_n \in \Val^{\Sp(n)}$ and
$r_{n+1,n}(\mu_{n+1})=\mu_n$ for $n=1,2,\ldots$. The inverse limit comes with
restriction maps 
\begin{align*}
 r_n:\Val^{\Sp(\infty)} & \to \Val^{\Sp(n)},\\
(\mu_1,\mu_2,\ldots) & \mapsto \mu_n
\end{align*}
which satisfy $r_n=r_{n+1,n} \circ r_{n+1}$ for all $n$.

Similarly, we denote by 
\begin{align*}
 \Val^{\Sp(\infty)\cdot \U(1)} & :=\lim_{\longleftarrow} \Val^{\Sp(n)\U(1)}\\
 \Val^{\Sp(\infty)\cdot \Sp(1)} & :=\lim_{\longleftarrow} \Val^{\Sp(n)\Sp(1)}
\end{align*}
the corresponding inverse limits. 

\proof[Proof of Theorem \ref{thm_local_global}]
We will prove the statement only in the case $G=\Sp(n)$, the two other cases are
similar. 

We have to show that the restriction map 
\begin{displaymath}
 r_n:\Val_k^{\Sp(\infty)}  \to  \Val_k^{\Sp(n)}
\end{displaymath}
is surjective for all $k$ and $n$ and injective for $n \geq k$. The case $k=0$
is trivial (each of the spaces $\Val_0^{\Sp(n)}$ and $\Val_0^{\Sp(\infty)}$ is
spanned by the Euler characteristic $\chi$). We thus assume that $k>0$. 

{\bf Injectivity:}
Let us first show that $r_n$ is injective for $n \geq k$. For this, it is enough
to show that $r_{n+1,n} :\Val_k^{\Sp(n+1)} \to \Val_k^{\Sp(n)}$ is injective for
all $n \geq k$.

Let $\mu \in \Val_k^{\Sp(n+1)}$ be a valuation whose restriction to
$\mathbb{H}^n$ vanishes. Take a space $E \in \Gr_k(\mathbb{H}^{n+1})$. Let
$E^\mathbb{H}$ be the quaternionic vector space generated by $E$. Note that
$\dim_\mathbb{H} E^\mathbb{H} \leq k \leq n$. 

We may find an element $g \in \Sp(n+1)$ with $g(E^\mathbb{H}) \subset
\mathbb{H}^n$. Then $F:=g(E) \in \Gr_k(\mathbb{H}^n)$ and hence $\mu|_F=0$.
Since $E$ and $F$ are in the same $\Sp(n+1)$-orbit and $\mu$ is
$\Sp(n+1)$invariant, it follows that $\mu|_E=0$. Since $E$ was arbitrary, the
Klain function of $\mu$ vanishes. By Theorem \ref{thm_klain}, $\mu$ vanishes.
This proves injectivity of $r_{n+1,n}$ and hence injectivity of $r_n$ for $n
\geq k$. 

{\bf Surjectivity:}

Let $n$ be arbitrary. We have the following maps
\begin{displaymath}
\xymatrixcolsep{2pc}
\xymatrix{\Omega^{k,k-1}(S\mathbb{H}^n)^{\Sp(n)} \ar_{*_1}^{\cong}[r] & 
\Omega^{4n-k,k-1}(S\mathbb{H}^n)^{\Sp(n)} \ar[r]^<<<<{\nc} & 
\Val_{4n-k}^{\Sp(n)} \ar[r]^{\cong}_{\F} & \Val_k^{\Sp(n)}}
\end{displaymath}
The first map is $*_1$, the Hodge star operator acting on the first component of
$\Omega^{*,*}(SV)^{\Sp(n)}=(\Lambda^*V \otimes \Omega^*(S^{4n-1}))^{\Sp(n)}$.

The second map is integration over the normal cycle from \eqref{eq_nc_onto}.
Since $k>0$, this map is surjective by \eqref{eq_inclusion_valg_smooth} and
\eqref{eq_nc_onto}. The third map is the Alesker-Fourier transform $\F$ (see
Theorem \ref{thm_alesker_fourier}), which is is an isomorphism. The composition
of these maps will be denoted by 
\begin{displaymath}
 q_n: \Omega^{k,k-1}(S\mathbb{H}^n)^{\Sp(n)} \twoheadrightarrow \Val_k^{\Sp(n)}, 
\end{displaymath}
this is a surjective map. 

Recall the definition of the surjective map $\tilde u_n$ from Proposition \ref{prop_free_algebra_up}. We claim that the following diagram commutes: 
\begin{displaymath}
\xymatrixcolsep{2pc}
\xymatrix{ & \Omega^{k,k-1}(S\mathbb{H}^{n+1})^{\Sp(n+1)} \ar[dd]^{\tilde
\iota_{n,n+1}^*} \ar@{->>}[r]^<<<<<{q_{n+1}} & \Val_k^{\Sp(n+1)} \ar[dd]^{r_{n+1,n}}\\
\Omega^{k,k-1}(S\mathbb{H}^\infty)^{\Sp(\infty)} \ar@{->>}[ur]^{\tilde u_{n+1}}
\ar@{->>}[dr]^{\tilde u_n}
& &\\
& \Omega^{k,k-1}(S\mathbb{H}^{n})^{\Sp(n)} \ar@{->>}[r]^<<<<<{q_n} & \Val_k^{\Sp(n)}
}
\end{displaymath}

The commutativity of the left hand triangle follows from the definitions. To
prove the commutativity of the right hand square, let $E \in
\Gr_k(\mathbb{H}^n)$, $\omega \in \Omega^{k,k-1}(S\mathbb{H}^{n+1})^{\Sp(n+1)}$
and let $B_E$ (resp. $B_{E^\perp}$, $S_E$) be the unit ball in $E$ (resp. in
$E^\perp$, resp. the unit sphere). Let us abbreviate $\iota:=\iota_{n,n+1},
\tilde \iota:=\tilde \iota_{n,n+1}, r:=r_{n+1,n}$. Then
\begin{align*}
\kl_{r \circ q_{n+1} \omega}(E) & =\kl_{q_{n+1}\omega}(\iota E)\\
& = \kl_{\nc(*_1\omega)}((\iota E)^\perp)\\
& = \frac{1}{\omega_{4n+4-k}} \int_{\nc(B_{(\iota E)^\perp})} *_1\omega\\
& = \frac{1}{\omega_{4n+4-k}} \int_{B_{(\iota E)^\perp} \times S_{(\iota E)}} *_1
\omega\quad \text{ since } \omega \text{ is of bi-degree } (k,k-1)\\
& = \frac{1}{\omega_k} \int_{B_{(\iota E)} \times S_{(\iota E)}} \omega. 
\end{align*}

On the other hand,
\begin{align*}
\kl_{q_n \circ \tilde \iota^* \omega}(E) & = \kl_{\nc(*_1 \tilde
\iota^*\omega)}(E^\perp)\\
& = \frac{1}{\omega_{4n-k}} \int_{\nc(B_E^\perp)} *_1 \tilde \iota^* \omega\\
& = \frac{1}{\omega_{4n-k}} \int_{B_E^\perp \times S_E} *_1 \tilde \iota^* \omega
\quad \text{ since } \tilde \iota^* \omega \text{ is of bi-degree } (k,k-1)\\
& = \frac{1}{\omega_k} \int_{B_E \times S_E} \tilde \iota^*\omega\\
& = \frac{1}{\omega_k} \int_{B_{(\iota E)} \times S_{(\iota E)}} \omega.  
\end{align*}
The claim now follows from Klain's embedding theorem \ref{thm_klain}.

The claim implies that for each $\omega \in
\Omega^{k,k-1}(S\mathbb{H}^{\infty})^{\Sp(\infty)}$, the element 
\begin{displaymath}
 \mu_\omega=(\mu_1,\mu_2,\ldots) \in \Val_k^{\Sp(\infty)}, \mu_n:=q_n \circ
\tilde u_n(\omega)
\end{displaymath}
is well-defined. 

Now we get surjectivity of $r_n$: let $\tau \in \Val_k^{\Sp(n)}$. Since
$\tilde u_n$ and
$q_n$ are surjective (see Proposition \ref{prop_free_algebra_up}), there exists
$\omega \in
\Omega^{k,k-1}(S\mathbb{H}^{\infty})^{\Sp(\infty)}$ with $\tau=q_n \circ
\tilde u_n(\omega)=r_n(\mu_\omega)$.  
\endproof

\proof[Proof of Theorem \ref{thm_main_asymptotic}]
The proof in the three cases $G=\Sp(\infty),
\Sp(\infty)\cdot \U(1), \Sp(\infty)\cdot \Sp(1)$ is similar and we first give an
outline. We let $G(n)$ be the corresponding group (i.e. if $G=\Sp(\infty)$
then $G(n)=\Sp(n)$ etc.). Set 
\begin{align*}
b_k & :=\dim (\Lambda^k \h^\infty \otimes \C)^G,\\
b_ {k,l} & :=\dim \Omega_h^{k,l}(S\h^\infty)^{\bar G}
\end{align*}
and 
\begin{align*}
b_k^n & :=\dim (\Lambda^k \h^n \otimes \C)^{G(n)},\\
b_ {k,l}^n & :=\dim \Omega_h^{k,l}(S\h^n)^{\bar G(n)}.
\end{align*}

We have the following symmetries (compare also Lemma \ref{lemma_char_rkl})
\begin{displaymath}
b^n_k=b^n_{4n-k}, \quad b^n_{k,l}=b^n_{l,k}=b^n_{4n-1-k,l}=b^n_{k,4n-1-l}. 
\end{displaymath}
Moreover, 
\begin{displaymath}
 b_{k,l}=\lim_{n \to \infty} b_{k,l}^n, \quad b_k=\lim_{n \to \infty} b_k^n.
\end{displaymath}
In fact, the first sequence stabilizes at $n=k+l+1$, the second at $n=k$ (see
Propositions \ref{prop_free_algebra_up} and \ref{prop_free_algebra_down}).
Let us introduce the formal power series 
\begin{align*}
 f^G(x) & =\sum_{k=0}^\infty b_k x^k,\\
 h^G(x,y) & =\sum_{k,l=0}^\infty b_{k,l} x^ky^l.
\end{align*}
It turns out that in all three cases there is a rational function $\tilde h^G$
such that
\begin{equation} \label{eq_splitting_rational}
 (x-y)h^G(x,y)=\tilde h^G(x,xy)-\tilde h^G(y,xy).
\end{equation}

Let us fix some $n > k+l$. By Theorems \ref{thm_local_global},
\ref{thm_alesker_fourier} and Corollary \ref{cor_exact_seq_group}, we have 
\begin{align*}
\dim \Val_k^{G} & = \dim \Val_k^{G(n)} \\
& = \dim \Val_{4n-k}^{G(n)} \\
& = \sum_{l=0}^{k-1} (-1)^{k+l+1}\dim  \Omega_p^{4n-k,l}(S\h^n)^{\bar
G(n)}+(-1)^k \dim (\Lambda^{4n-k} \h^n \otimes \C)^{G(n)}\\
& = \sum_{l=0}^{k-1}
(-1)^{k+l+1}(b^n_{4n-k,l}-b^n_{4n-k-1,l-1})+(-1)^k b^n_{4n-k}\\
& = \sum_{l=0}^{k-1} (-1)^{k+l+1} (b^n_{k-1,l}-b^n_{k,l-1}) + (-1)^k b^n_k\\
& = \sum_{l=0}^{k-1} (-1)^{k+l+1} (b_{k-1,l}-b_{k,l-1}) + (-1)^k b_k
\end{align*}
We multiply by $x^k$ and sum over $k$ to obtain
\begin{align*}
\sum_{k=0}^\infty \dim \Val_k^G x^k&
= - \left(\sum_{l=0}^{k-1} (-1)^l (b_{k-1,l}-b_{k,l-1})\right)(-x)^k +\sum_{k=0}^\infty (-1)^k b_kx^k\\
& = -\tilde h^G(-x,x)+f^G(-x).
\end{align*}

Let us compute the function $f^G$. By Proposition
\ref{prop_free_algebra_down} and \eqref{eq_symmetric_powers}, $(\Lambda^k
\h^\infty \otimes \C)^{\Sp(\infty)}=0$ if $k$ is odd and   
\begin{displaymath}
 (\Lambda^k \h^\infty \otimes \C)^{\Sp(\infty)} \simeq \Sym^{k/2}V_2
\simeq \bigoplus_{l=0}^{\lfloor k/4\rfloor} V_{k-4l}
\end{displaymath}
if $k$ is even. Therefore 
\begin{align*}
 \dim (\Lambda^k \h^\infty \otimes \C)^{\Sp(\infty)} & =\sum_{l=0}^{\left\lfloor
\frac{k}{4}\right\rfloor}(k-4l+1)=\binom{\frac{k}{2}+2}{2}\\
\dim (\Lambda^k \h^\infty \otimes \C)^{\Sp(\infty)\cdot \U(1)} & =\left\lfloor
\frac{k}{4}\right\rfloor+1\\
\dim (\Lambda^k \h^\infty \otimes \C)^{\Sp(\infty)\cdot \Sp(1)} & = \left\{
\begin{array}{c
c} 1 & k \equiv 0 \mod 4\\0 & k \not\equiv 0 \mod 4\end{array}\right.,
\end{align*}
from which we deduce that 
\begin{align*}
f^{\Sp(\infty)}(x) &
=\frac{1}{(1-x^2)^3} \\
f^{\Sp(\infty) \cdot \U(1)}(x) &
=\frac{1}{(1-x^2)(1-x^4)} \\
f^{\Sp(\infty) \cdot \Sp(1)}(x) &
=\frac{1}{1-x^4}. 
\end{align*}

Next, we compute the function $h^G$. By Proposition \ref{prop_free_algebra_up},
there is an isomorphism of
$\SU(2)$-representation
\begin{multline} \label{eq_sl2_decomposition_forms}
 \Omega_h^{k,l}(S\h^\infty)^{\overline{\Sp(\infty)}}=\bigoplus
\Sym^{k_1}V_2 \otimes \Sym^{k_2}V_2
\otimes \\
\Sym^{k_3}V_0 \otimes \Sym^{k_4}V_2 \otimes
\Sym^{k_5}V_2 \otimes \Sym^{k_6}V_2
\end{multline}
where the sum is over all tuples $(k_1,\ldots,k_6)$ with $k_1+k_3+k_5+2k_6=k$
and $k_2+k_3+2k_4+k_5=l$. 

{\bf Claim:} The $\SU(2)$-representation $\Omega_h^{k,l}(S\h^\infty)^{\overline{\Sp(\infty)}}$ does not contain $V_{2s+1}$ for any $s$ and contains $V_{2s}$ with a multiplicity which is given by the coefficient of $x^ky^l$ in the power series
\begin{multline} \label{eq_power_series}
\frac{1}{(x-y)^2(x+y)(1-x^4)(1-y^4)(1-x^2y^2)(1-xy)} \cdot \\
\bigg( \frac{(1+x)^2(1+y)(1+x^3)(1+y^2)(x^2+y)(x^2y+1)}{(1-xy)(1-x^3y)} x^{2s}\\
 -\frac{(1+x)^2(1+y)^2(1+y^2)(1+x^2)(x^2y+1)(y^2x+1)(x+y)}{(-1+xy^3)(-1+x^3y)}
x^sy^s\\
 +\frac{(1+y)^2(1+y^3)(1+x^2)(y^2+x)(y^2x+1)(1+x)}{(1-xy)(1-xy^3)} y^{2s}\bigg).
\end{multline}

We defer the technical proof of the claim, which only uses the
Clebsch-Gordan rule, to the appendix. Multiplying \eqref{eq_power_series} by
$2s+1$ and summing over all $s$ we obtain 
\begin{displaymath}
h^{\Sp(\infty)}(x,y)= \frac{1}{(1-x)^3(1-xy)^4(1-y)^3}
\end{displaymath}

Setting
\begin{displaymath}
 \tilde h^{\Sp(\infty)}(x,y):=\frac{x(1-xy)}{(1-x)^3(1-y)^7},
\end{displaymath}
\eqref{eq_splitting_rational} is satisfied. It follows that 
\begin{align*}
 \sum_{k=0}^\infty \dim \Val_k^{\Sp(\infty)}x^k & =f^{\Sp(\infty)}(-x)-\tilde
h^{\Sp(\infty)}(-x,x)\\
& =\frac{x^4-3x^3+6x^2-3x+1}{(1-x)^7(1+x)^3}.
\end{align*}

Similarly, summing \eqref{eq_power_series} over all $s$, we obtain  
\begin{displaymath}
 h^{\Sp(\infty)\cdot \U(1)}(x,y)=\frac{A}{B},
\end{displaymath}
where
\begin{align*}
A & := -2y^2+y-1-x(3y-1)(y+1)(y^2-y+1)+x^2(y^4-5y^3-2)\\
& \quad -x^3y(2y^4+5y-1)+x^4y(y-3)(1+y)(y^2-y+1)-x^5y^3(y^2-y+2)\\
B & := (-1+y)^2(-1+xy^3)(-1+x)^2(-1+xy)^3(xy+1)(1+y^2)(1+x^2)(-1+x^3y)
\end{align*}

Now one can check that 
\begin{displaymath}
 \tilde h^{\Sp(\infty)\cdot \U(1)}(x,y):=\frac{\tilde A(x,y)}{\tilde B(x,y)}
\end{displaymath}
with 
\begin{align*}
\tilde A(x,y) & := -x(1+3y^2+3y^3+y^5)+x^2(3y^6-3y+y^2-y^4-y^5+1)\\
& \quad +x^3(-4y^2+2y-2y^6+2y^7-2+4y^5)\\
& \quad +x^4(-y^5+y^2-3y+3y^6+y^3-y^7)+x^5(y^7+y^2+3y^5+3y^4)\\
\tilde B(x,y) & := (-1+y)^6(-1+x^2y)(-1+x)^2(1+x^2)(1+y)^2(1+y^2)(y^2+y+1)
\end{align*}
satisfies \eqref{eq_splitting_rational}. 

It follows that 
\begin{align*}
 \sum_{k=0}^\infty \dim \Val_k^{\Sp(\infty)\cdot \U(1)}x^k & =f^{\Sp(\infty)\cdot
\U(1)}(-x)-\tilde h^{\Sp(\infty)\cdot
\U(1)}(-x,x)\\
& =\frac{x^6-2x^5+2x^4+2x^2-2x+1}{(x^2+1)(x^2+x+1)(1+x)^2(1-x)^6}.
\end{align*}

Finally, setting $s=0$ in \eqref{eq_power_series} we obtain   
\begin{displaymath}
 h^{\Sp(\infty)\cdot \Sp(1)}(x,y)=\frac{A}{B},
\end{displaymath}
where
\begin{align*}
A & :=(-1+x^2y^2-3x^3y^3-x^4y^4+x^6y^6-x-x^2-xy-y-y^2+x^6y^4+x^3+y^3\\
& \quad
+2x^3y-x^3y^5+2xy^3-xy^5-2x^2y^4-x^5y^3+x^4y^6-y^4x^3-x^4y^3-2x^4y^2-x^5y\\
& \quad +2x^3y^2+2x^2y+2y^2x+x^4y+2x^2y^3+y^4x)(x-y)\\
B & := (-1+xy^3)(-1+x^3y)(-1+x^2y^2)(y^3-y^2+y-1)(x^3-x^2+x-1)(-1+xy)^2
\end{align*}

Now one can check that 
\begin{displaymath}
 \tilde h^{\Sp(\infty)\cdot \Sp(1)}(x,y):=\frac{\tilde A(x,y)}{\tilde B(x,y)}
\end{displaymath}
with 
\begin{align*}
\tilde A(x,y) & := x(1+y)(y^5+y^4+2y^3+2y^2+1)+x^2(2y^5+y^3-1)\\
& \quad +x^3(-y^7-y^5+y^4+y^3+y^2-1)+x^4(y+1)(2y^4+3y^2+y+1)\\
& \quad -x^5y(-1+y)(1+y)(y^2+1)(y^2+y+1)\\
\tilde B(x,y) & := (-x^3+x^2-x+1)(-1+y^2)(-1+y)^2\\
& \quad \cdot (-1+x^2y)(y^6+y^5+y^4-y^2-y-1)
\end{align*}
satisfies \eqref{eq_splitting_rational}. 

It follows that 
\begin{align*}
 \sum_{k=0}^\infty \dim \Val_k^{\Sp(\infty)\cdot \Sp(1)}x^k & =f^{\Sp(\infty)\cdot
\Sp(1)}(-x)-\tilde h^{\Sp(\infty)\cdot \Sp(1)}(-x,x)\\
& =\frac{1+x^3+2x^4+x^5}{(x^2+1)(x^2+x+1)(1+x)^2(1-x)^4}.
\end{align*}
\endproof

\begin{Corollary}
The dimension of the space of $k$-homogeneous translation $G$-invariant
continuous valuations is given by 
 \begin{align} 
\dim \Val_k^{\Sp(\infty)} & \label{eq_explicit_dimension_spn}
=\frac{(k+4)(k+2)(2k^4+24k^3+100k^2+168k+405+315(-1)^k)}{5760}\\
\dim \Val_k^{\Sp(\infty)\cdot \U(1)} & \label{eq_explicit_dimension_spnu1}
=\frac{1}{8}\left(i^k+(-i)^k\right)+\frac{4\sqrt{3}(-1)^{k+1}}{81}
\sin\left(\frac{\pi k}{3}\right) \nonumber \\
& \quad +\frac{(15+5k)(-1)^k}{64}+\frac{33}{64}+\frac{421}{960}k+\frac{5}{24}k^2
\nonumber\\
& \quad +\frac{7}{108}k^3+\frac{1}{96}k^4+\frac{1}{1440}k^5\\
 \dim \Val_k^{\Sp(\infty)\cdot \Sp(1)} & \label{eq_explicit_dimension_spnsp1}
=\frac{3}{16}\left(i^k+(-i)^k\right)+\frac{2\sqrt{3}(-1)^{k+1}}{27}
\sin\left(\frac{\pi k}{3}\right) \nonumber \\
& \quad
+\frac{(5+k)(-1)^k}{32}+\frac{15}{32}+\frac{37}{96}k+\frac{3}{16}k^2+\frac{5}{
144}k^3
\end{align}
\end{Corollary}

\proof
This follows from Theorem \ref{thm_main_asymptotic} by a standard technique. 
\endproof

Here we list the first few values of these dimensions. 
\begin{equation} \label{eq_asymptotic_dim_examples}
 \begin{array}{c | c | c | c}
  k & \dim \Val_k^{\Sp(\infty)} & \dim \Val_k^{\Sp(\infty)\cdot \U(1)} & \dim
\Val_k^{\Sp(\infty)\cdot \Sp(1)} \\ \hline
0 & 1 & 1 & 1\\
1 & 1 & 1 & 1\\
2 & 7 & 3 & 2\\
3 & 14 & 6 & 4\\
4 & 42 & 14 & 8\\
5 & 84 & 24 & 11\\ 
6 & 182 & 44 & 17\\
7 & 330 & 72 & 24\\
8 & 603 & 117 & 34 \\
9 & 1001 & 177 & 44\\
10 & 1645 & 265 & 58
 \end{array}
\end{equation}

\section{A character formula for $\Val^{\spn}$}
\label{sec_dimformula}

The aim of this section is the proof of Theorem \ref{thm_main_thm}. Let us fix
some notation first. Let $V$ be a quaternionic vector space of dimension $n$,
endowed with a quaternionic Hermitian form $K$. The group $\Sp(V,K)$ is then
isomorphic to $\Sp(n)$. Next, we fix some unit vector $v_0 \in V$ and write 
\begin{displaymath}
 V=\R \cdot v \oplus U \oplus \tilde V, 
\end{displaymath}
where $\tilde V$ is the quaternionic complement of $v_0 \cdot \h$ and $U$ is
three-dimensional. 

The stabilizer of $\Sp(n) \cong \Sp(V,K)$
at $v_0$ is $\Sp(n-1) \cong \Sp(\tilde V,K|_{\tilde V})$, it acts trivially on
$U$. As
$\SU(2)$-representation $U \simeq V_2$, the adjoint representation.

In the following, we will consider real and complex vector spaces. In order to
distinguish between both cases, we will write $\Lambda_\R$ and $\Lambda_\C$
referring to the base field $\R$ and $\C$ respectively. 

\begin{Lemma} \label{lemma_dim_invariant_forms}
 \begin{displaymath}
 \dim (\Lambda_\C^{k_1}V \otimes \ldots \otimes
\Lambda_\C^{k_s}V)^{\Sp_{2n}\mathbb{C}}=\sum_{\lambda} K_{\lambda
\mu},
\end{displaymath}
where $\mu=(k_1,\ldots,k_s)$ and the sum is over all even Young diagrams
$\lambda$  with $\lambda_1 \leq 2n$. 
\end{Lemma}

\proof
We can consider $V$ as a representation of the larger group
$\GL(2n,\C)$. By \eqref{eq_dec_prod_alt}, we have 
\begin{displaymath}
 \Lambda_\C^{k_1}V \otimes \ldots \otimes \Lambda_\C^{k_s}V=\bigoplus_{\tilde
\nu}
K_{\tilde \nu \mu} \Gamma_\nu,
\end{displaymath}
where $\nu$ ranges over all partitions of depth $\leq 2n$, $\tilde \nu$ is the
conjugate partition of $\nu$ and $K_{\tilde \nu \mu}$ is the Kostka number. 

Proposition \ref{prop_spherical} and the fact that $K_{\lambda \mu}=0$ unless
$\lambda \trianglerighteq \mu$ thus imply that 
\begin{align*}
 \dim\left(\Lambda_\C^{k_1}V \otimes \ldots \otimes
\Lambda_\C^{k_s}V\right)^{\Sp_{2n}\mathbb{C}} &
=\sum_{\substack{\tilde
\nu \text{ even},\\ \text{depth } \nu \leq 2n}} K_{\tilde \nu
\mu}=\sum_{\substack{\lambda \text{ even
}\\ \lambda_1 \leq 2n}} K_{\lambda \mu}.
\end{align*}
\endproof

\begin{Lemma} \label{lemma_char_down}
Let 
\begin{displaymath}
 S_k:=(\Lambda_\R^k V)^{\spn} \otimes \C,
\end{displaymath}
which is a $\SU(2)$-representation. Then 
\begin{displaymath} 
 \sum_{k=0}^{4n} \ch(S_k) x^k=E_n(x),
\end{displaymath}
where $E_n(x)$ is the function defined in \eqref{eq_def_en}.
\end{Lemma}

\proof
We have seen in Section \ref{sec_global} that 
\begin{displaymath}
 (\Lambda^k V)^{\spn} \otimes \C=\bigoplus_{k_1+k_2=k} (\Lambda_\C^{k_1}V
\otimes \Lambda_\C^{k_2}V^*)^{\Sp_{2n}\C}
\end{displaymath}
Elements of $(\Lambda_\C^{k_1}V \otimes
\Lambda_\C^{k_2}V^*)^{\Sp_{2n}\C}$ are of weight $k_1-k_2$, hence 
\begin{align*}
\sum_{k=0}^{4n} \ch(S_k) x^k & = \sum_{k=0}^{4n} \ch((\Lambda_\R^kV)^{\spn}
\ \otimes \C) x^k \\
& =\sum_{k=0}^{4n} \sum_{k_1+k_2=k} \dim (\Lambda_\C^{k_1}V
\otimes
\Lambda_\C^{k_2}V^*)^{\Sp_{2n}\C} t^{k_1-k_2}x^k\\
& = \sum_{\lambda} \sum_{k_1,k_2} K_{\lambda,(k_1,k_2)} (xt)^{k_1}(xt^{-1})^{k_2} \\
& = \sum_\lambda s_\lambda(xt,xt^{-1})\\
& = E_n(x),
\end{align*}
where $\lambda$ in the above sums ranges over all even Young diagrams of
depth $\leq 2$ with $\lambda_1 \leq 2n$. 
\endproof

Note that $S_k=0$ if $k$ is odd and that $S_k=S_{4n-k}$. 

\begin{Lemma} \label{lemma_char_rkl}
Set 
\begin{displaymath}
 R_{k,l} :=\left(\Lambda_\R^{k,l}(\tilde V \oplus \tilde V)\right)^{\Sp(n-1)} 
\otimes \C,
\end{displaymath}
which is an $\SU(2)$-representation. 
Then for each $m$ 
\begin{displaymath}
\sum_{k=0}^{2m} \ch R_{k,2m-k}x^k= F_{n-1,m}(x).
\end{displaymath}
\end{Lemma}

\proof
\begin{align*}
\sum_{k=0}^{2m} \ch(R_{k,2m-k})x^k & = \sum_{k=0}^{2m}
\sum_{\substack{k_1+k_2=k\\ l_1+l_2=2m-k}} 
\dim \left(\Lambda_\C^{k_1}\tilde V \otimes \Lambda_\C^{k_2}\tilde V^* \otimes
\Lambda_\C^{l_1}\tilde V
\otimes \Lambda_\C^{l_2}\tilde V^*\right)^{\Sp_{2n-2}\C} \cdot \\
& \quad \cdot x^kt^{k_1-k_2+l_1-l_2}\\
& = \sum_\lambda \sum_k
\sum_{\substack{k_1+k_2=k\\ l_1+l_2=2m-k}} 
{K_{\lambda,(k_1,k_2,l_1,l_2)}}x^k t^{k_1-k_2+l_1-l_2}\\
& = \sum_{\substack{\lambda}}
\sum_{\substack{k_1,k_2,l_1,l_2\\ k_1+k_2+l_1+l_2=2m}} 
{K_{\lambda,(k_1,k_2,l_1,l_2)}}(xt)^{k_1}(xt^{-1})^{k_2} t^{l_1-l_2}\\
& = \sum_{\lambda}
s_\lambda(tx,t^{-1}x,t,t^{-1})\\
& = F_{n-1,m}(x).
\end{align*}
Here $\lambda$ ranges over all even Young diagrams of depth $\leq 4$, weight
$2m$ with $\lambda_1 \leq 2n-2$.
\endproof

Note that we have isomorphisms of
$\SU(2)$-representations $R_{k,l} = R_{l,k}=R_{k,4n-4-l}$ and that $R_{k,l}=0$
if $k+l$ is odd. 

\begin{Proposition} \label{prop_sl2_dec_valk}
In $R\SU(2)$, we have the following equation 
\begin{align} \label{eq_sl2_dec_valk}
\Val_k^G & =S_k -R_{k,k-2} -R_{k,k-4} +  R_{k-1,k-1} \nonumber \\
& \quad + (-V_4+V_2) R_{k-1,k-3} \nonumber \\
& \quad + (V_4-V_2+V_0) R_{k-2,k-2} \nonumber \\
& \quad -R_{k-2,k-4} + R_{k-3,k-3}.
\end{align}
\end{Proposition}

\proof
By Proposition \ref{prop_exact_seq}, 
\begin{equation} \label{eq_alt_sum_valk}
\Val_k^{\spn} = (-1)^k (\Lambda_\R^kV)^{\spn} \otimes \C
+\sum_{l=0}^{4n-k-1}(-1)^{k+l+1}\Omega^{k,l}_p(SV)^{\spn}.
\end{equation}

By definition,
\begin{align*}
\Omega^{k,l}_p(SV)^{\spn} &
=\Omega^{k,l}_h(SV)^{\spn}-\Omega^{k-1,l-1}_h(SV)^{\spn}\\
& = \left[\Lambda_\R^k (U \oplus \tilde V)^* \otimes \Lambda_\R^l(U \oplus
\tilde V)^* \otimes \C\right]^{\Sp(n-1)}\\
& - \left[\Lambda_\R^{k-1} (U \oplus \tilde V)^* \otimes \Lambda_\R^{l-1}(U
\oplus
\tilde V)^* \otimes \C\right]^{\Sp(n-1)}. 
\end{align*}

Note that for all $j$
\begin{displaymath}
\Lambda_\R^j (U \oplus \tilde V)=\Lambda_\R^j \tilde V \oplus V_2 \otimes
\Lambda_\R^{j-1}\tilde V \oplus V_2 \otimes \Lambda_\R^{j-2}\tilde V \oplus
\Lambda_\R^{j-3}\tilde V.
\end{displaymath}

In the alternating sum \eqref{eq_alt_sum_valk}, most of the terms cancel
out and 
we get, for any given $k$, 
 equation \eqref{eq_sl2_dec_valk}.
\endproof

\proof[Proof of Theorem \ref{thm_main_thm}]

Using Proposition \ref{prop_sl2_dec_valk} and Lemmas \ref{lemma_char_down} and
\ref{lemma_char_rkl} we obtain
\begin{align*}
\sum_{k=0}^{4n} \ch(\Val_k^G)x^k & = \sum_{k=0}^{4n} \ch(S_k)x^k -\sum_{k=0}^{4n} \ch(R_{k,4n-k-2})x^k -\sum_{k=0}^{4n} \ch(R_{k,4n-k})x^k \\
& \quad +\sum_{k=0}^{4n} \ch(R_{k-1,4n-k-3})x^k + \sum_{k=0}^{4n} \ch(-V_4+V_2) \ch(R_{k-1,4n-k-1}) x^k \\
& \quad + \sum_{k=0}^{4n} \ch(V_4-V_2+V_0) \ch(R_{k-2,4n-k-2})x^k-\sum_{k=0}^{4n} \ch(R_{k-2,4n-k})\\
& \quad + \sum_{k=0}^{4n} \ch(R_{k-3,4n-k-1})x^k\\
& = E_n(x) -F_{n-1,2n}(x)-(1+x(t^4+t^{-4})+x^2)F_{n-1,2n-1}(x)\\
& \quad +x(1+ x(t^4+1+t^{-4})+ x^2)F_{n-1,2n-2}(x)  
\end{align*}
\endproof

As an example, we work out the cases $n=1$ and $n=2$. 

\subsubsection*{Case $n=1$}
The only Young diagram with $\lambda_1\leq 2n-2=0$ is $\lambda=(0,\ldots,0)$ and
it has weight $0$. Hence $F_{0,m}=0$ for $m \neq 0$ and $F_{0,0}=1$. The
polynomial $E_1$ equals 
\begin{displaymath}
E_1(x)=x^4+x^2(t^2+1+t^{-2})+1 
\end{displaymath}
From our formula it follows that 
\begin{align*}
 \sum_{k=0}^4 \ch(\Val_k^{\Sp(1)})x^k & =x^4+x^2(t^2+1+t^{-2})+1+x(1+ x(t^4+1+t^{-4})+ x^2)\\
 & = x^4+x^3+x^2(t^4+t^2+2+t^{-2}+t^{-4})+x+1,
\end{align*}
hence  
\begin{align*}
\Val_k^{\Sp(1)} & = V_0 \quad \text{ for } k =0,1,3,4\\
\Val_2^{\Sp(1)} & = V_4+V_0
\end{align*}

These decompositions have been obtained in \cite{bernig_sun09} (note that $\Sp(1)=\SU(2)$).

\subsubsection*{Case $n=2$}

The polynomials $E_2$ and $F_{1,m}, m=2,3,4$ are given as follows:
\begin{align*}
E_2 & = x^8+x^6(t^2+1+t^{-2})+x^4(t^4+t^2+1+t^{-2}+t^4)+x^2(t^2+1+t^{-2})+1\\
F_{1,2} & = x^4+x^3(t^2+2+t^{-2})+x^2(t^4+2t^2+4+2t^{-2}+t^{-4})+x(t^2+2+t^{-2})+1\\
F_{1,3} & = x^4(t^{-2}+1+t^2)+x^3(t^{-2}+2+t^2)+x^2(t^{-2}+1+t^2) \\
F_{1,4} & = x^4. 
\end{align*}

Putting these values in the formula from Theorem \ref{thm_main_thm}, we obtain 

\begin{align*}
 \sum_{k=0}^8 \ch(\Val_k^{\Sp(2)})x^k & =x^8+x^7+x^6(t^4+t^2+3+t^{-2}+t^{-4})\\
 & \quad +x^5(2t^4+2t^2+5+2t^{-2}+2t^{-4})\\
 & \quad +x^4(t^8+t^6+4t^4+4t^2+9+4t^{-2}+4t^{-4}+t^{-6}+t^{-8})\\
 & \quad +x^3(2t^4+2t^2+5+2t^{-2}+2t^{-4})\\
 & \quad +x^2(t^4+t^2+3+t^{-2}+t^{-4})+x+1,
\end{align*}
which implies that 
\begin{align}
\Val_0^{\Sp(2)} & = \Val_1^{\Sp(2)} = \Val_7^{\Sp(2)} = \Val_8^{\Sp(2)} =
V_0 \nonumber \\
\Val_2^{\Sp(2)} & = \Val_6^{\Sp(2)}= V_4+2V_0 \nonumber\\
\Val_3^{\Sp(2)} & = \Val_5^{\Sp(2)} = 2V_4+3V_0 \nonumber\\
\Val_4^{\Sp(2)} & = V_8+3V_4+5V_0. \label{eq_dec_val_sp2}
\end{align}

Similar computations yield the decomposition of each $\Val_k^{\Sp(n)}$, from which we may compute the following dimensions:
 
\begin{displaymath}
 \begin{array}{c | c} n & \dim \Val_k^{\spn}, k=0,\ldots,4n\\ \hline
1 & 1,1,6,1,1\\
2 & 1,1,7,13,29,13,7,1,1\\
3 & 1,1,7,14,41,71,111,71,41,14,7,1,1\\
4 & 1,1,7,14,42,83,169,259,344,259,169,83,42,14,7,1,1\\
5 & 1,1,7,14,42,84,181,317,532,742,903,742,532,317,181,84,42,14,7,1,1
 \end{array}
\end{displaymath}

\begin{displaymath}
 \begin{array}{c | c} n & \dim \Val_k^{\spnu}, k=0,\ldots,4n\\ \hline
1 & 1,1,2,1,1\\
2 & 1,1,3,5,9,5,3,1,1\\
3 & 1,1,3,6,13,19,25,19,13,6,3,1,1\\
4 & 1,1,3,6,14,23,39,53,64,53,39,23,14,6,3,1,1\\
5 & 1,1,3,6,14,24,43,67,98,124,141,124,98,67,43,24,14,6,3,1,1
 \end{array}
\end{displaymath}

\begin{displaymath}
 \begin{array}{c | c} n & \dim \Val_k^{\spnsp}, k=0,\ldots,4n\\ \hline
1 & 1,1,1,1,1\\
2 & 1,1,2,3,5,3,2,1,1\\
3 & 1,1,2,4,7,8,9,8,7,4,2,1,1\\
4 & 1,1,2,4,8,10,14,16,18,16,14,10,8,4,2,1,1\\
5 & 1,1,2,4,8,11,16,21,26,28,30,28,26,21,16,11,8,4,2,1,1
 \end{array}
\end{displaymath}

Note that, for $k \leq n$, these numbers are consistent
with \eqref{eq_asymptotic_dim_examples} and Theorem
\ref{thm_local_global}. 

\begin{appendix}
\section{Proof of \eqref{eq_power_series}}
\label{sec_power_series}

Let $R:=R\SU(2)$ and define  
\begin{align*}
 T_s & := \sum_{k=0}^\infty V_0 (xy)^k \in R[[x,y]]\\
 T_0 & :=\sum_{k=0}^\infty \Sym^kV_2 y^{2k} \in R[[y]]\\
  T_1 & :=\sum_{k=0}^\infty \Sym^kV_2 (xy)^k \in R[[x,y]]\\
 T_2 & :=\sum_{k=0}^\infty \Sym^kV_2 x^{2k} \in R[[x]]\\
B & := V_0+V_2x+V_2x^2+V_0x^3 \in R[[x]]\\
G & := V_0+V_2y+V_2y^2+V_0y^3 \in R[[y]].
\end{align*}

Then \eqref{eq_sl2_decomposition_forms} can be rewritten as 
\begin{equation} \label{eq_sl2_decomposition_forms_power_series}
 \sum_{k,l} \Omega^{k,l}_h(S\h^\infty)^{\overline{\Sp(\infty)}}x^ky^l \cong B \cdot G \cdot T_s \cdot T_0 \cdot T_1 \cdot T_2.
\end{equation}

Set 
\begin{displaymath}
 p :=\sum_{k=0}^\infty V_{2k}x^{2k}, q :=\sum_{k=0}^\infty V_{2k}y^{2k},  
 r :=\sum_{k=0}^\infty V_{2k}(xy)^k \in R[[x,y]].
\end{displaymath}
 
\begin{Lemma} \label{lemma_multiplication_sl2_ring}
 \begin{align}
 T_0 & = \frac{q}{1-y^4} \label{eq_t0_term}\\
T_1 & =\frac{r}{1-x^2y^2} \label{eq_t1_term}\\
T_2 &= \frac{p}{1-x^4} \label{eq_t2_term} \\
T_s & = \frac{V_0}{1-xy} \label{eq_ts_term}\\
(x^2-y^2) p \otimes q & = \frac{x^2(1+y^2)p - y^2(1+x^2)q}{1-x^2y^2} \label{eq_pq_term}\\
(x-y) p \otimes r & = \frac{x(1+xy)p - y(1+x^2)r}{1-x^3y} \label{eq_pr_term} \\
(x-y) r \otimes q & = \frac{x(1+y^2)r - y(1+xy)q}{1-xy^3} \label{eq_rq_term} \\
x B \otimes p & = (1+x)^2(1+x^3)p-(1+x)(1+x^2)V_0 \label{eq_Bp_term}\\ 
 y G \otimes q & = (1+y)^2(1+y^3)q-(1+y)(1+y^2)V_0. \label{eq_Gq_term}
\end{align}
\end{Lemma}

\proof
We have 
\begin{displaymath}
 \Sym^kV_2=\sum_{s=0}^{\lfloor \frac{k}{2} \rfloor} V_{2k-4s}
\end{displaymath}
and hence 
\begin{multline*}
 T_0  = \sum_{k=0}^\infty \Sym^kV_2 y^{2k} = \sum_{k=0}^\infty \sum_{s=0}
^{\lfloor \frac{k}{2} \rfloor} V_{2k-4s} y^{2k}= \sum_{s=0}^\infty y^{4s}
\sum_{k=2s}^\infty V_{2k-4s}y^{2k-4s}\\
= \frac{1}{1-y^4} \sum_{k=0}^\infty V_{2k}y^{2k}= \frac{q}{1-y^4}.
\end{multline*}

This shows \eqref{eq_t0_term}; the proofs of \eqref{eq_t1_term} and
\eqref{eq_t2_term} are similar. Equation \eqref{eq_ts_term} is trivial. 

Next, for \eqref{eq_pq_term} we compute 
\begin{align*}
 (x^2-y^2) p \otimes q & = (x^2-y^2) \sum_{k,l=0}^\infty V_{2k} \otimes V_{2l} x^{2k}y^{2l}\\
& =  (x^2-y^2) \sum_{k,l=0}^\infty \sum_{s=|k-l|}^{k+l} V_{2s} x^{2k}y^{2l}\\
& = (x^2-y^2) \sum_{s=0}^\infty V_{2s} \sum_{k,l:|k-l|\leq s \leq k+l} x^{2k}y^{2l}\\
& = (x^2-y^2) \sum_{s=0}^\infty V_{2s} \sum_{k=0}^\infty x^{2k} \sum_{l=|s-k|}^{k+s}y^{2l}\\
& = (x^2-y^2) \sum_{s=0}^\infty V_{2s} \sum_{k=0}^\infty x^{2k} \frac{y^{2|k-s|}-y^{2(k+s+1)}}{1-y^2}\\
& = \frac{x^2-y^2}{1-y^2}\sum_{s=0}^\infty V_{2s} \left(\frac{y^{2s+2}-x^{2s+2}}{y^2-x^2}+\frac{x^{2s+2}y^2}{1-x^2y^2}-\frac{y^{2s+2}}{1-x^2y^2}\right)\\
& =\frac{x^2(y^2+1)}{(1-x^2y^2)}p - \frac{y^2(1+x^2)}{(1-x^2y^2)}q\\
\end{align*}
Equations \eqref{eq_pr_term} and \eqref{eq_rq_term} are proved in a similar way. 

Finally, Equations \eqref{eq_Bp_term} and \eqref{eq_Gq_term} follow easily from the Gordan-Clebsch rule. 
\endproof

From the lemma, we deduce that 
\begin{align*}
B \cdot G \cdot  T_s \cdot T_0 \cdot T_1 \cdot T_2 & = c\frac{(1+x)(1+x^3)(1+y^2)(x^2+y)(x^2y+1)}{(1-xy)(x^2-y^2)(1-x^3y)(x-y)} p\\
& \quad -c\frac{(1+y^2)(1+x^2)(1+x)(1+y)(x^2y+1)(y^2x+1)}{(-1+xy^3)(y-x)^2(-1+x^3y)}r\\
& \quad +c\frac{(1+y)(1+y^3)(1+x^2)(y^2+x)(y^2x+1)}{(1-xy)(x^2-y^2)(1-xy^3)(x-y)}q
\end{align*}
where 
\begin{displaymath}
 c:=\frac{(1+x)(1+y)}{(1-xy)(1-x^4)(1-y^4)(1-x^2y^2)}.
\end{displaymath}
Comparing the coefficient of $V_{2s}$ on both sides, taking into account \eqref{eq_sl2_decomposition_forms_power_series}, yields \eqref{eq_power_series}.

\end{appendix}


\end{document}